\documentclass[a4paper,12pt]{article}
\usepackage{mathptmx}      
\usepackage[babel=true]{csquotes}
\usepackage{amsmath,amsfonts,graphicx,amssymb} 
\usepackage{float}
\usepackage{array}
\usepackage{hyperref}
\usepackage{framed}
\usepackage[normalem]{ulem}
\usepackage{algorithmic}
\usepackage{arydshln} 

\DeclareMathAlphabet{\mathcal}{OMS}{cmsy}{m}{n}

\usepackage{algorithm}

\newcommand{\Vct}[1]{{\boldsymbol{\mathbf{#1}}}}
\newcommand{\Mat}[1]{{\boldsymbol{\mathbf{#1}}}}

\newcommand{\PtThreeD}[1]{\mathsf{\uppercase{#1}}}
\newcommand{\PtTwoD}[1]{\mathsf{\textbf{\lowercase{#1}}}}
\newcommand{\Transpose}[1]{#1^\top}

\newcommand{\Bydef}{\stackrel{\vartriangle}{=}}

\newcommand{\R}{{\mathbb{R}}}
\newcommand{\N}{{\mathbb{N}}}
\newcommand{\K}{{\mathit{K}}}

\newcommand{\Nn}{\mathbb{N}^n}

\newcommand{\RNn}{\mathbb{R}^{\mathbb{N}^n}}

\newcommand{\y}{\mathbf{y}}
\newcommand{\fmin}{f^{\star}}
\newcommand{\fmom}{\widehat{f}}
\newcommand{\fmomt}{\widehat{f}_t}

\newcommand{\borK}{\mathcal{M}_+(K)}

\begin{document}

\title{Rank-Constrained Fundamental Matrix Estimation by Polynomial Global Optimization Versus the Eight-Point Algorithm}

\author{Florian~Bugarin,$^{1}$ Adrien~Bartoli,$^{2}$
 	Didier~Henrion,$^{3,4}$\\ Jean-Bernard~Lasserre,$^{3,5}$ Jean-Jos\'e~Orteu,$^{1}$ Thierry~Sentenac,$^{1,3}$
}

\footnotetext[1]{Universit\'e de Toulouse; INSA, UPS, Mines Albi, ISAE; ICA (Institut Cl\'ement Ader); Campus Jarlard, F-81013 Albi, France}

\footnotetext[2]{Universit\'é d'Auvergne, ISIT, UMR 6284 CNRS / Universit\'e d'Auvergne, Clermont-Ferrand, France}

\footnotetext[3]{LAAS-CNRS, Universit\'e de Toulouse, Toulouse, France}

\footnotetext[4]{Faculty of Electrical Engineering, Czech Technical University in Prague}

\footnotetext[5]{Institut de Math\'ematiques de Toulouse, Universit\'e de Toulouse, Toulouse, France}

\maketitle

\begin{abstract}
The fundamental matrix can be estimated from point matches.
The current gold standard is to bootstrap the eight-point algorithm and two-view projective bundle adjustment.
The eight-point algorithm first computes a simple linear least squares solution by minimizing an algebraic cost and then projects the result to the closest rank-deficient matrix. 
We propose a single-step method that solves both steps of the eight-point algorithm.
Using recent results from polynomial global optimization, our method finds the rank-deficient matrix that exactly minimizes the algebraic cost. In this special case, the optimization method is reduced to the resolution of very short sequences of convex linear problems which are computationally efficient and numerically stable. 
The current gold standard is known to be extremely effective but is nonetheless outperformed by our rank-constrained method for bootstrapping bundle adjustment.
This is here demonstrated on simulated and standard real datasets.
With our initialization, bundle adjustment consistently finds a better local minimum (achieves a lower reprojection error) and takes less iterations to converge.
\\


\end{abstract}


\section{Introduction}
The fundamental matrix has received a great interest in the computer vision community (see for instance~\cite{longuet-higgins:81a,salvi:03,tsai:84a,hartley:95a,Torr02,Chojnacki03,BartoliSturm:04}). 
This $(3\times3)$ rank-two matrix encapsulates the epipolar geometry, the projective motion between two uncalibrated perspective cameras, and serves as a basis for 3D reconstruction, motion segmentation and camera self-calibration, to name a few. 
Given $n$ point matches $(\PtTwoD{q}_i,\PtTwoD{q}_i')$, $i=1,\dots,n$ between two images, the fundamental matrix may be estimated in two phases.
The initialization phase finds some suboptimal estimate while the refinement phase iteratively minimizes an optimal but nonlinear and nonconvex criterion.
The gold standard uses the {\em eight-point algorithm} and {\em projective bundle adjustment} for these two phases, respectively.
A `good enough' initialization is necessary to avoid local minima at the refinement phase as much as possible.
The main goal of this article is to improve the current state of the art regarding the initialization phase. We here focus on input point matches that do not contain mismatches (a pair of points incorrectly associated). The problem of mismatches has been specifically addressed by the use of robust methods in the literature.

The eight-point algorithm follows two steps~\cite{longuet-higgins:81a}.
In its first step, it relaxes the rank-deficiency constraint and solves the following convex problem: 
\begin{equation}
  \tilde{\Mat{F}} = \underset{\Mat{F} \in \mathbb{R}^{3\times3}}{\arg\min}\; C(\Mat{F}) \mbox{ s.t. } \lVert \Mat{F} \rVert^2 =1,
  \label{equ:alg_dist_estimation}
\end{equation}
where $C$ is a convex, linear least squares cost, hereinafter called the {\em algebraic cost}:
\begin{equation}
C(\Mat{F}) = \displaystyle \sum_{i=1}^n \left( \Transpose{\PtTwoD{q}_i'}\Mat{F}\PtTwoD{q}_i \right) ^2.
\end{equation}
This minimization is subject to the normalization constraint $\lVert \Mat{F} \rVert^2 =1$.
This is to avoid the trivial solution $\Mat{F}=\Mat{0}$.
Normalization will be further discussed in section~\ref{sec:gloptipoly}.
The estimated matrix $\tilde{\Mat{F}}$ is thus {\em not} a fundamental matrix yet.
In its second step, the eight-point algorithm computes the closest rank-deficient matrix to $\tilde{\Mat{F}}$ as:
\begin{equation}
  \Mat{F}_{\text{8pt}} = \underset{\Mat{F} \in \mathbb{R}^{3\times3}}{\arg\min}\; \displaystyle \lVert \Mat{F} - \tilde{\Mat{F}} \rVert^2 \mbox{ s.t. } \det(\Mat{F})=0.
  \label{equ:rank_correction}
\end{equation}
Both steps can be easily solved.
The first step is a simple linear least squares problem and the second step is solved by nullifying the least singular value of $\tilde{\Mat{F}}$.
It has been shown~\cite{hartley:95a} that this simple algorithm performs extremely well in practice, provided that the image point coordinates are {\em standardized} by simply rescaling them so that they lie in $[-\sqrt{2};\sqrt{2}]^2$.

Our main contribution in this paper is an approach that solves for the fundamental matrix minimizing the algebraic cost.
In other words, we find the global minimum of:
\begin{equation}
  \tilde{\Mat{F}} = \underset{\Mat{F} \in \mathbb{R}^{3\times3}}{\arg\min}\; C(\Mat{F}) \mbox{ s.t. } \det(\Mat{F})=0 \mbox{ and } \lVert \Mat{F} \rVert^2 =1.
  \label{equ:cons_alg_dist_estimation}
\end{equation}
Perhaps more importantly, we also quantify the impact that each of $\Mat{F}_{\text{8pt}}$ and $\tilde{\Mat{F}}$ has when used as an initial estimate in Bundle Adjustment. Each initial estimate will lead Bundle Adjustment to its own refined estimate. The two final estimates may thus be different since, as the difference between the two initial estimates grows larger, the probability that they lie in different basins of attraction increases. Our measure quantifies:
\begin{enumerate}
 \item how far are these two basins of attractions,
 \item how many iterations will Bundle Adjustment take to converge.
\end{enumerate}
The proposed algorithm uses polynomial global optimization~\cite{Lasserre:01,gloptipoly:02}. Previous attempts~\cite{Xiao:10,Chesi:02,KahlHenrion:07} in the literature differ in terms of optimization strategy and parameterization of the fundamental matrix. None solves problem~(\ref{equ:cons_alg_dist_estimation}) optimally for a general parameterization: they either do not guarantee global optimality~\cite{Chesi:02,KahlHenrion:07} or prescribe some camera configurations~\cite{Xiao:10,Chesi:02,KahlHenrion:07} (requiring typically that the epipole in the first camera does not lie at infinity). Furthermore, the main criticism made to the optimization method we use is the resolution of a hierarchy of convex linear problems of increasing size, which is computationally ineffective and numerically unstable. The proposed solution overcomes this drawback: experiments show that, in most of cases, the proposed algorithm only requires solving the second relaxation of the sequences used in the exposed optimization method.

Our experimental evaluation on simulated and real datasets compares the difference between the eight-point algorithm and ours used as initialization to bundle adjustment.
We observe that {\em (i)} bundle adjustment consistently converges within less iterations with our initialization and {\em (ii)} bundle adjustment always achieves an equal or lower reprojection error with our initialization.
We provide numerous examples of real image pairs from standard datasets.
They all illustrate practical cases for which our initialization method allows bundle adjustment to reach a better local minimum than the eight-point algorithm.
 
\section{State of the Art}
Accurately and automatically estimating the fundamental matrix from a pair of images is a major research topic.
We first review a four-class categorization of existing methods, and specifically investigate the details of existing global methods.
We finally state the improvements brought by our proposed global method.

\subsection{Categorizing Methods}
A classification of the different methods in three categories --linear, iterative and robust-- was proposed in~\cite{salvi:03}.
Linear methods directly optimize a linear least squares cost.
They include the eight-point algorithm~\cite{longuet-higgins:81a}, SVD resolution~\cite{salvi:03} and variants~\cite{tsai:84a,hartley:95a,Torr02,Chojnacki03}. 
Iterative methods iteratively optimize a nonlinear and nonconvex cost.
They require, and are sensitive to the quality of, an initial estimate.
The first group of iterative methods minimizes the distances between points and epipolar lines~\cite{salviPhd:99,PeiChen:10}.
The second group minimizes some approximation of the reprojection error~\cite{luong:96a,Zhang:98,Hengel02anew,Chojnacki02anew}. 
The third group of methods minimizes the reprojection error, and are equivalent to two-view projective bundle adjustment.
Iterative methods typically use a nonlinear parameterization of the fundamental matrix which guarantees that the rank-deficiency constraint is met. 
For instance, a minimal 7-parameter update can be used over a consistent orthogonal representation~\cite{BartoliSturm:04}.
Finally, robust methods estimate the fundamental matrix while classifying each point match as inlier or outlier.
Robust methods use M-Estimators~\cite{Hartley-Zisserman:03}, LMedS (median least squares)~\cite{Zhang:98} or RANSAC (random sampling consensus)~\cite{Torr97thedevelopment}.
Both LMedS and RANSAC are stochastic.

To these three categories, we propose to add a fourth one: {\em global methods}.
Global methods attempt at finding the global minimum of a nonconvex problem.
Convex relaxations have been used to combine a convex cost with the rank-deficiency constraint~\cite{Chesi:02}. 
However, these relaxations do not converge to a global minimum and the solution's optimality is not certified. 

\subsection{Global Methods}
In theory, for a constrained optimization problem, global optimization methods do not require an initial guess and may be guaranteed to reach the global minimum, thereby certifying optimality. Such global methods can be separated in two classes. The methods of the first class describe the search space as exhaustively as possible in order to test as many candidate solutions as possible. Following this way, there are methods such as Monte-Carlo sampling, which test random elements satisfisying constraints, and reactive tabu search~\cite{Glover:86,Jaumard:90}, which continues searching even after a local minimum has been found. The major drawback of these methods is mainly in the prohibitive computation time required to have a sufficiently high probability of success. Moreover, even in case of convergence, there is no certificate of global optimality. Contrary to the methods of the first class, methods lying in the second class provide a certificate of global optimality using the mathematical theory from which they are built. Branch and Bound algorithms~\cite{Lang:60} or global optimization by interval analysis~\cite{hansen:03,Kearfott:09} are some examples. However, although these methods can be faster than those of the first category, their major drawback is their lack of generality. Indeed, these methods are usually dedicated to one particular type of cost function because they use highly specific computing mechanisms to be as efficient as possible. A review of global methods may be found in~\cite{Weise:07}.  

A good deal of research has been conducted over the last few decades on applying global optimization methods in order to solve polynomial minimization problems under polynomial constraints. The major drawback of these applications has been the difficulty to take constraints into account. But, by solving simplified problems, these approaches have mainly been used to find a starting point for local iterative methods. However, recent results in the areas of convex and polynomial optimization have facilitated the emergence of new approaches. These have attracted great interest in the computer vision community. In particular, global polynomial optimization~\cite{Lasserre:01,gloptipoly:02} has been used in combination with a finite-epipole nonlinear parameterization of the fundamental matrix~\cite{Xiao:10}. This method does not consequently cover camera setups where the epipole lies at infinity. A global convex relaxation scheme~\cite{Lasserre:01,gloptipoly:02} was used to minimize the Sampson distance~\cite{KahlHenrion:07}. Because this implies minimizing a sum of many rational functions, the generic optimization method had to be specifically adapted and lost the property of certified global optimality.

\subsection{The Proposed Method}
The proposed method lies in the fourth category: it is a global method.
Similarly to the eight-point algorithm, it minimizes the algebraic cost, but explicitly enforces the nonlinear rank-deficiency constraint. Contrarily to previous global methods~\cite{Xiao:10,Chesi:02,KahlHenrion:07}, the proposed method handles all possible camera configurations (it does not make an assumption on the epipoles being finite or infinite) and certifies global optimality. Moreover, the presented algorithm is based on the resolution of a very short sequence of convex linear problems and is therefore computationally efficient.

A large number of attemps to introduce global optimization have been made in the literature. In~\cite{Chesi:02}, a dedicated hierarchy of convex relaxations is defined in order to globally solve the problem of fundamental matrix estimation. The singularity constraint is taken into account by introducing additional optimization variables, not introduced directly in the problem description. The resulting optimization algorithm is not generic and, contrary to Lasserre's hierarchy, there is no proof that the sequence of solutions of this specific hierarchy converges to the global minimum.

In~\cite{Xiao:10}, Lasserre's hierarchy is used jointly with the introduction of the singularity constraint in the problem description. However, the normalization constraint $\|\Mat{F}\|^2=1$ is replaced by fixing, a priori, one of the $\Mat{F}$ coefficients to 1. Hence the other $\Mat{F}$ coefficients are not bounded. Consequently, there is no guarantee that the sequence of solutions $(\fmomt)_{t\in \mathbb{N}}$ converges to the global minimum.

In~\cite{kahl:05} the authors minimize the Sampson distance by solving a specific hierarchy of convex relaxations built upon an epigraph formulation. This specific hierarchy is obtained by extending Lasserre's hierarchy, where the contraints are linear (LMIs), to a hierarchy of convex relaxations where constraints are polynomial (polynomial matrix inequality - PMI). Consequently, no asymptotic convergence of the hierarchy to the global optimum can be guaranteed. Moreover, even if the rank correction is achieved using the singularity constraint, the normalization is replaced by setting, a priori, the coefficient $f_{33}$ to 1 and so the other coefficients of $\Mat{F}$ are not bounded.

Finally, in a very recent work~\cite{Zheng:13}, the algebraic error is globally minimized thanks to the resolution of seven subproblems. Each subproblem is reduced to a polynomial equation system solved via a Gr\"{o}bner basis solver. The singularity constraint is satisfied thanks to the right epipole parametrization. Although this parametrization ensures that $\Mat{F}$ is singular while using the minimum number of parameters, this method is not practical since it would be necessary to solve 126 subproblems in order to cover all the 18 possible parameter sets~\cite{Zhang:98}. Therefore it is preferable to introduce the singularity constraint directly in the problem description rather than via some parametrization of $\Mat{F}$.
\\


\section{Polynomial Global Optimization}\label{sec:gloptipoly}
\subsection{Introduction}
Given a real-valued polynomial $f(x) : \mathbb{R}^n \rightarrow \mathbb{R}$, we are interested in solving the problem:
\begin{eqnarray}\label{eq:prob_opti_globale}
\fmin = \underset{x \in \K}{\inf}\; f(x)
\end{eqnarray}
where $\K\subseteq \mathbb{R}^n$ is a (not necessarily convex) compact set defined by polynomial inequalities:  $g_j(x)\geq 0, \;\; j = 1,\ldots,m$. Our optimization method is based on an idea first described in \cite{lasserre00}. It consists in reformulating the nonconvex global optimization problem~\eqref{eq:prob_opti_globale} as the equivalent convex linear programming problem:
 \begin{eqnarray}\label{eq:prob_opti_globale_mesure}
 \fmom = \underset{\mu \in \mathcal{P}(\K)}{\inf}\;
\displaystyle \int_\K f(x) d\mu,
\end{eqnarray}
where $\mathcal{P}(\K)$ is the set of probability measures supported on $\K$. Note that this reformulation is true for any continuous function (not necessarily polynomial) and any compact set $K\subseteq \mathbb{R}^n$. Indeed, as $f^\star \leq f(x)$, then $f^\star \leq \int_\K f d\mu $ and thus $f^\star  \leq  \widehat{f}$. Conversely, if $x^\star$ is a global minimizer of~\eqref{eq:prob_opti_globale}, then the probability measure $\mu^\star \Bydef \delta_{x^\star}$ (the Dirac at $x^\star$) is admissible for~\eqref{eq:prob_opti_globale_mesure}. Moreover, because $\fmom$ is a solution of~\eqref{eq:prob_opti_globale_mesure}, the following inequality holds: $\int_\K f(x) d\mu \geq \fmom, \; \forall \mu \in \mathcal{P}(\K)$ and thus $f^\star=\int_{K} f(x) \, \delta_{x^\star} \geq \widehat{f}$. Instead of optimizing over the finite-dimensional euclidean space $\K$, we optimize over the infinite-dimensional set of probability measures $\mathcal{P}(\K)$. Thus, Problem~\eqref{eq:prob_opti_globale_mesure} is, in general, not easier to solve than Problem~\eqref{eq:prob_opti_globale}. However, in the special case of $f$ being a polynomial and $\K$ being defined by polynomial inequalities, we will show how Problem~\eqref{eq:prob_opti_globale_mesure} can be reduced to solving an (generically finite) sequence of convex linear matrix inequality (LMI) problems.

\subsection{Notations and Definitions}
First, given vectors $\alpha = (\alpha_1,\ldots,\alpha_n)^\top \in \mathbb{N}^n$  and $x = (x_1,\ldots,x_n)^\top \in \mathbb{R}^n$, we define the \textit{monomial} $x^\alpha$ by:  
\begin{eqnarray}
  x^\alpha &\Bydef& x_1^{\alpha_1}x_2^{\alpha_2}\ldots x_n^{\alpha_n}
\end{eqnarray}
and its \textit{degree} by $\mathrm{deg}(x^\alpha)\Bydef \displaystyle \lVert \alpha \rVert_1 = \sum_{i=1}^n \alpha_i$. For $t\in \N$, we define $\mathbb{N}^n_t$ the space of the $n$-dimensional integer vector with a norm lower than $t$ as:
\begin{eqnarray}
  \mathbb{N}^n_t &\Bydef & \left\lbrace  \alpha \in \mathbb{N}^n \;\; | \;\;  \lVert \alpha \rVert_1 \leq t  \right\rbrace. 
\end{eqnarray}
Then, consider the family:
\begin{eqnarray}
\left\lbrace x^\alpha\right\rbrace _{\alpha \in \mathbb{N}_t^n} & = &\left\lbrace 1,x_1,x_2,\ldots,x_n,x_1^2,x_1x_2,\ldots,x_1x_n,x_2x_3,\ldots,x_n^2,\ldots,x_1^t,\ldots,x_n^t\right\rbrace 
\end{eqnarray}
of all the monomials $x^\alpha$ of degree at most $t$, which has dimension $s(t)\Bydef \dfrac{(n+t)!}{t!n!}$. Those monomials form the canonical basis of the vector space $\mathbb{R}_t[x]$ of real-valued multivariate polynomials of degree at most $t$. Then, a polynomial $p \in \mathbb{R}_t[x]$ is understood as a linear combination of monomials of degree at most $t$:
\begin{equation}
p(x) = \sum_{\alpha \in \mathbb{N}_t^n} p_\alpha x^\alpha,
\end{equation}
and $\mathbf{p}\Bydef(p_\alpha)_{\lVert \alpha \rVert_1 \leq t}\in \R^{\mathbb{N}_t^n}\simeq \R^{s(t)}$ is the vector of its coefficients in the monomial basis $\left\lbrace x^\alpha\right\rbrace _{\alpha \in \mathbb{N}_t^n}$. Its degree is equal to $\mathrm{deg}(p)\Bydef \max{\left\lbrace \lVert \alpha \rVert_1 \; | \; p_\alpha\neq0\right\rbrace } $ and $\mathrm{d}_p$ denotes the smallest integer not lower than $\dfrac{\mathrm{deg}(p)}{2}$.
\\

\noindent \textit{Example}:
The polynomial
\begin{eqnarray}\label{example1}
    x\in \R^2 &\mapsto& p(x)=1+2x_2 +3x^2_1 +4x_1x_2
\end{eqnarray}
has a vector of coefficients $p\in \R^6$ with entries $p_{00}=1$, $p_{10}=0$, $p_{01}=2$, $p_{20}=3$, $p_{11}=4$ and $p_{02}=0$.
\\

Next, given $\mathbf{y}=(y_\alpha)_{\alpha \in \N^n} \in \R^{\N^n}$, we define the \textit{Riesz functional} $L_\mathbf{y}$ by the linear form: 
\begin{equation}
 L_\mathbf{y} :  \begin{array}[t]{lll}
  \R\left[x \right] & \rightarrow & \mathbb{R} \\ 
 		   p=\displaystyle\sum_{\alpha \in \mathbb{N}^n} p_\alpha x^\alpha  & \rightarrow & \mathbf{y}^\top \mathbf{p}= \displaystyle \sum_{\alpha \in \mathbb{N}^n} p_\alpha y_\alpha.
  \end{array}
\end{equation}
Thus, the Riesz functional can be seen as an operator that linearizes polynomials.
\\

\noindent \textit{Example:}
For the polynomial~\eqref{example1}, the Riesz functional reads
\begin{eqnarray} \label{example2}
     p(x)=1+2x_2 +3x^2_1 +4x_1x_2 &\mapsto& L_\mathbf{y}(p) = y_{00}+2y_{01}+3y_{20}+4y_{11}.
\end{eqnarray}
\\
 
For $t\in \N$ and $\y \in \R^{\N^n_{2t}}$, the matrix $M_t(\mathbf{y})$ of size $s(t)$ defined by:
\begin{eqnarray}
		(M_t(\mathbf{y}))_{\alpha,\beta}= L_\mathbf{y}(x^\alpha x^\beta )= y_{\alpha+\beta} \quad \forall \alpha,\beta \in \mathbb{N}_{t}^n 
\end{eqnarray}
is called \textit{the moment matrix of order} $t$ \textit{of} $\y$. By construction, this matrix is symmetric and linear in $\y$. Then, given $q\in\R_t\left[x \right]$ and $\mathbf{q}\in \R^{\N^n_{t}}$ the vector of its coefficients in the monomial basis, the vector: 
\begin{eqnarray}
		q\mathbf{y}&\Bydef&M_t(\mathbf{y})\mathbf{q} \;\in\; \R^{\N^n_{t}}
\end{eqnarray}
is called \textit{the shifted vector} with respect to $\mathbf{q}$. $M_t(q\mathbf{y})$, the moment matrix of order $t$ of $q\y$, is called \textit{the localizing matrix of degree} $t$ \textit{of} $q$. This matrix is also symmetric and linear in $\y$.
\\

\noindent\textit{Example}:
If $n=2$ then: 
\begin{equation} \label{example3}
\begin{array}{l}
 M_0(\mathbf{y})=y_{00},\\
M_1(\mathbf{y})=\left(
			\begin{array}{c:cc}
				y_{00}  \;\;&\;\; y_{10} & y_{01} \\ \hdashline 
				y_{10} \;\; &\;\; y_{20} & y_{11} \\ 
				y_{01}  \;\;&\;\; y_{11} & y_{02} \\
			\end{array}
		\right),\\
M_2(\mathbf{y})=\left(
			\begin{array}{c:cc:ccc}
				y_{00}  \;\; & \;\;y_{10} & y_{01} \;\;&\;\; y_{20}  & y_{11} & y_{02} \\ \hdashline 
				y_{10}  \;\; & \;\;y_{20} & y_{11} \;\;&\;\; y_{30}  & y_{21} & y_{12} \\
				y_{01}  \;\; & \;\;y_{11} & y_{02} \;\;&\;\; y_{21}  & y_{12} & y_{03} \\ \hdashline 
				y_{20}  \;\; & \;\;y_{30} & y_{21} \;\;&\;\; y_{40}  & y_{31} & y_{22} \\
				y_{11}  \;\; & \;\;y_{21} & y_{12} \;\;&\;\; y_{31}  & y_{22} & y_{13} \\
				y_{02}  \;\; & \;\;y_{12} & y_{03} \;\;&\;\; y_{22}  & y_{13} & y_{04} \\
			\end{array}
		\right),	
\end{array}	
\end{equation}
and if $q(x)=a+2x_1^2+3x_2^2$ then:
\begin{eqnarray} \label{example4}
M_1(q\mathbf{y})=\left(
			\begin{array}{cc:ccc}
				ay_{00}+2y_{20}+3y_{02}&  & \;\;ay_{10}+2y_{30}+3y_{12}& & ay_{01}+2y_{21}+3y_{03} \\ \hdashline 
				ay_{10}+2y_{30}+3y_{12}&  & \;\;ay_{20}+2y_{40}+3y_{22}& & ay_{11}+2y_{31}+3y_{13} \\ 
				ay_{01}+2y_{21}+3y_{03}&  & \;\;ay_{11}+2y_{31}+3y_{13}& & ay_{02}+2y_{22}+3y_{04} \\
			\end{array}
		\right).	
\end{eqnarray}
\\

Finally, recall that a symmetric matrix $F\in \mathbb{S}^n$ is positive semidefinite, denoted by $F\succeq0$, if and only if  $x^\top F x\geq0, \; \forall x \in \R^n$ or equivalently, if and only if the minimum eigenvalue of $F$ is non-negative. A \textit{linear matrix inequality} (\textit{LMI}) is a convex constraint: 
\begin{eqnarray}
		F_0+\displaystyle \sum_{k=1}^n x_kF_k&\succeq &0,
\end{eqnarray}
on a vector $x\in \R^n$, where matrices $F_k \in \mathbb{S}^m$, $k = 0,\ldots,n$ are given.

\subsection{Optimization Method}
Let $f$ be a real-valued multivariate polynomial, Problem~\eqref{eq:prob_opti_globale_mesure} can be reduced to a convex linear programming problem. Indeed, if $f(x)=\sum_{\alpha \in \N^n} f_\alpha x^\alpha$ then:
\begin{equation}\label{eg:int_L}
  \int_K f\; d\mu = \int_K \sum_{\alpha \in \mathbb{N}^n} f_\alpha x^\alpha d\mu = \sum_{\alpha \in \mathbb{N}^n} f_\alpha \int_K x^\alpha d\mu = L_\mathbf{y}(f)  
\end{equation}
where each coordinate $y_\alpha$ of the infinite sequence $\y \in \R^{\mathbb{N}^n}$ is equal to $\displaystyle \int_K x^\alpha \, \mu(dx)$, also called \textit{the moment of order} $\alpha$. Consequently, if $f$ is polynomial, then Problem~\eqref{eq:prob_opti_globale_mesure} is equivalent to:  
\begin{eqnarray}\label{eq:pmom2}
\fmom&=&\begin{array}[t]{l}
	\inf  \; L_\mathbf{y}(f)\\
	\text{s.t. }\begin{array}[t]{l}
			y_0=1 \\
			\mathbf{y} \in \mathcal{M}_\textit{K}.
			\end{array}
\end{array}
\end{eqnarray}
with:
\begin{eqnarray}
\mathcal{M}_\textit{K}&\Bydef&
\left\lbrace \mathbf{y} \in \RNn \;\mid \; \exists \mu \in \borK \text{ such that }  y_\alpha = \int_K x^\alpha d\mu\; \; \forall \alpha \in \Nn \right\rbrace,
\end{eqnarray}
and $\borK$ is the space of finite Borel measures supported on $\K$. Remark that the constraint $y_0=1$ is added in order to impose that if $\mathbf{y} \in \mathcal{M}_\textit{K}$ then $\mathbf{y}$ represents a measure in $\mathcal{P}(\K)$ (and no longer in $\borK$). Although Problem~\eqref{eq:pmom2} is a convex linear programming problem, it is difficult to describe the convex cone $\mathcal{M}_\textit{K}$ with simple constraints on $\y$. But, the problem \enquote{$\y \in \mathcal{M}_\textit{K}$}, also called $\K$-\textit{moment problem}, is solved when $\K$ is a \textit{basic semi-algebraic set}, namely:
\begin{eqnarray}\label{def:semi-algebrique}
 K&\Bydef&\left\{ x \in \mathbb{R}^n \hspace{0.2cm}\vert\hspace{0.2cm} g_1(x)\geqslant0,\ldots,g_m(x)\geqslant0  \right\}
\end{eqnarray}
where $g_j \in \R[x], \; \forall j=1,\ldots m$. Note that $\K$ is assumed to be compact. Then, without loss of generality, we assume that one of the polynomial inequalities $g_j(x)\geqslant0$ is of the form $R^2-\lVert x \rVert_2^2\geqslant0$ where $R$ is a sufficiently large positive constant. In this particular case, $\mathcal{M}_\textit{K}$ can be modelled using LMI conditions on the matrices $M_t(\y)$ and $M_t(g_j\y), \;j=1,\ldots m$. More precisely, thanks to the Putinar Theorem~\cite{LasserreBook:10,Laurent:08}, we have: 
\begin{eqnarray}\label{lemme:Put}
\mathcal{M}_\textit{K} &=& \mathcal{M}_\succeq(g_1,\ldots,g_m),
\end{eqnarray}
where: 
\begin{equation}
\begin{array}{ccc}
 \mathcal{M}_\succeq(g_1,\ldots,g_m)&\Bydef &\left\lbrace \mathbf{y} \in \mathbb{R}^{\mathbb{N}^n} \vert \,
M_t(\mathbf{y}) \succeq 0,\, M_t(g_j\mathbf{y}) \succeq 0 \; \; \forall j =1,\ldotp \ldotp,m\; \; \forall t \in \N
\right\rbrace.
\end{array}
\end{equation}
Then, Problem~\eqref{eq:prob_opti_globale_mesure} is equivalent to: 
\begin{eqnarray}\label{eq:pmom3}
\fmom=\begin{array}[t]{l}
	\underset{\y \, \in \, \mathbb{R}^{\mathbb{N}^n}}{\inf}\displaystyle L_\y(f)\\
	\text{s.t. } \begin{array}[t]{l}
			y_0=1 \\
			M_t(\y) \succeq 0\\
			M_t(g_j\y) \succeq 0 \;\;\; j =1,\ldotp \ldotp,m \;\;\; \forall t \in \N.
			\end{array}
\end{array}
\end{eqnarray}
To summarize, if $f$ is polynomial and $\K$ a semi-algebraic set, then Problem~\eqref{eq:prob_opti_globale} is equivalent to a convex linear programming problem with an infinite number of linear constraints on an infinite number of decision variables. Now, for $t\geq \mathrm{d}_{\K}\Bydef\max(\mathrm{d}_f,\mathrm{d}_{g_1},\ldots,\mathrm{d}_{g_m})$ consider the finite-dimensional truncations of Problem~\eqref{eq:pmom3}:
\begin{equation}\label{def:fmomt}
 \mathcal{Q}_t\Bydef\left\lbrace \begin{array}{lll}\fmomt & \stackrel{\vartriangle}{=}&
				\begin{array}[t]{l}
					\underset{\y \, \in \, \mathbb{R}^{\mathbb{N}_{2t}^n}}{\min}\displaystyle L_\y(f) \\
					\text{ s.t. } \begin{array}[t]{l}
								y_0 =1 \\
								M_t(\y) \succeq 0,\\
					 			M_{t-\mathrm{d}_{g_j}}(g_j\y) \succeq 0\; \; \forall j \in  \left\lbrace 1,\ldots,m\right\rbrace.
							\end{array}
   				\end{array}
   			\end{array}
\right. 
\end{equation}
By construction, $\mathcal{Q}_t$, $t \in \mathbb{N}$ generates a hierarchy of LMI relaxations of Problem~\eqref{eq:pmom3}~\cite{Lasserre:01}, where each $\mathcal{Q}_t$, $t \in \mathbb{N}$, is concerned with moment and localizing matrices of fixed size $t$. Each relaxation~\eqref{def:fmomt} can be solved by using public-domain implementations of primal-dual interior point algorithms for semidefinite programming (SDP)~\cite{Borchers99csdp,dsdp-user-guide,Fujisawa95sdpa,sdpt3UserGuide,S98guide}. When \textit{the relaxation order} $t\in \N$ tends to infinity, we obtain the following results~\cite{Lasserre:01,Nie2012}: 
\begin{equation}
	\widehat{f}_{t} \leq \widehat{f}_{t+1} \leq  \widehat{f} \text{ and } \underset{t\rightarrow +\infty}{\lim} \widehat{f}_t  = \widehat{f} .
\end{equation} 
Practice reveals that this convergence is fast and very often finite, i.e. there exists a finite $t_0$ such that $\widehat{f}_t  = \widehat{f}$, $\forall t\geq t_0$. In fact, finite convergence is guaranteed in a number of cases (e.g. discrete optimization) and very recent results by Nie \cite{Nie2012} show that finite convergence is even generic as well as an optimal solution $\y^\star_t$ of~\eqref{def:fmomt}.
\\

\noindent \textit{Example}:
Consider the polynomial optimization problem
\begin{equation}
\begin{array}{lll}\fmom & =&
	\begin{array}[t]{ll}
	\displaystyle\min_{x \in \R^2} & -x_2 \\
	\mathrm{s.t.} & 3-2x_2-x_1^2-x_2^2 \geq 0 \\
	& -x_1-x_2-x_1x_2 \geq 0 \\
	& 1+x_1x_2 \geq 0.
	\end{array}
\end{array}	
\end{equation}
The first LMI relaxation $\mathcal{Q}_1$ is 
\begin{equation}
\begin{array}{lll}\widehat{f}_1 & =&
	\begin{array}[t]{ll}
	\displaystyle\min_{\y \in \R^6} & -y_{01} \\
	\mathrm{s.t.} & y_{00}=1\\
		      &\begin{pmatrix}
		         y_{00} & y_{10} & y_{01} \\
		         y_{10} & y_{20} & y_{11} \\
		         y_{01} & y_{11} & y_{02} 
		       \end{pmatrix}\succeq 0\\
		      & 3y_{00}-2y_{01}-y_{20}-y_{02} \geq 0 \\
		      & -y_{10}-y_{01}-y_{11} \geq 0 \\
		      & y_{00}+y_{11} \geq 0,
	\end{array}
\end{array}	
\end{equation}
and the second LMI relaxation $\mathcal{Q}_2$ is
\begin{equation}
\begin{array}{lll}\widehat{f}_2 & =&
	\begin{array}[t]{ll}
	\displaystyle\min_{\y \in \R^{15}} & -y_{01} \\
	\mathrm{s.t.} & y_{00}=1\\
		      &\begin{pmatrix}
			y_{00} & y_{10} & y_{01} & y_{20} & y_{11} & y_{02} \\ 
			y_{10} & y_{20} & y_{11} & y_{30} & y_{21} & y_{12} \\
			y_{01} & y_{11} & y_{02} & y_{21} & y_{12} & y_{03} \\  
			y_{20} & y_{30} & y_{21} & y_{40} & y_{31} & y_{22} \\
			y_{11} & y_{21} & y_{12} & y_{31} & y_{22} & y_{13} \\
			y_{02} & y_{12} & y_{03} & y_{22} & y_{13} & y_{04}
		       \end{pmatrix}\succeq 0,\\
		      &\begin{pmatrix}
				3y_{00}-2y_{01}-y_{20}-y_{02} & 3y_{10}-2y_{11}-y_{30}-y_{12} & 3y_{01}-2y_{02}-y_{21}-y_{03} \\ 
				3y_{10}-2y_{11}-y_{30}-y_{12} & 3y_{20}-2y_{21}-y_{40}-y_{22} & 3y_{11}-2y_{12}-y_{31}-y_{13} \\ 
				3y_{01}-2y_{02}-y_{21}-y_{03} & 3y_{11}-2y_{12}-y_{31}-y_{13} & 3y_{02}-2y_{03}-y_{22}-y_{04} 
			\end{pmatrix}\succeq 0\\
		      &\begin{pmatrix}
				-y_{10}-y_{01}-y_{11} & -y_{20}-y_{11}-y_{21} & -y_{11}-y_{02}-y_{12} \\ 
				-y_{20}-y_{11}-y_{21} & -y_{30}-y_{21}-y_{31} & -y_{21}-y_{31}-y_{21} \\ 
				-y_{11}-y_{02}-y_{12} & -y_{21}-y_{12}-y_{22} & -y_{12}-y_{03}-y_{13} 
			\end{pmatrix}\succeq 0\\
		      &\begin{pmatrix}
				y_{00}+y_{11} & y_{10}+y_{21} & y_{01}+y_{12} \\ 
				y_{10}+y_{21} & y_{20}+y_{31} & y_{11}+y_{22} \\ 
				y_{01}+y_{12} & y_{11}+y_{22} & y_{02}+y_{13} 
		      \end{pmatrix}\succeq 0.
	\end{array}
\end{array}	
\end{equation}
It can be checked that $\widehat{f}_1=-2\leq\widehat{f}_2=\fmom=-\frac{1+\sqrt{5}}{2}$. Note that the constraint $3-2x_2-x_1^2-x_2^2 \geq 0$ certifies boundedness of the feasibility set.
\\

\begin{figure*}[!htp]\label{fig:gloptipoly}
  \centering
   \includegraphics[height = 5.5in]{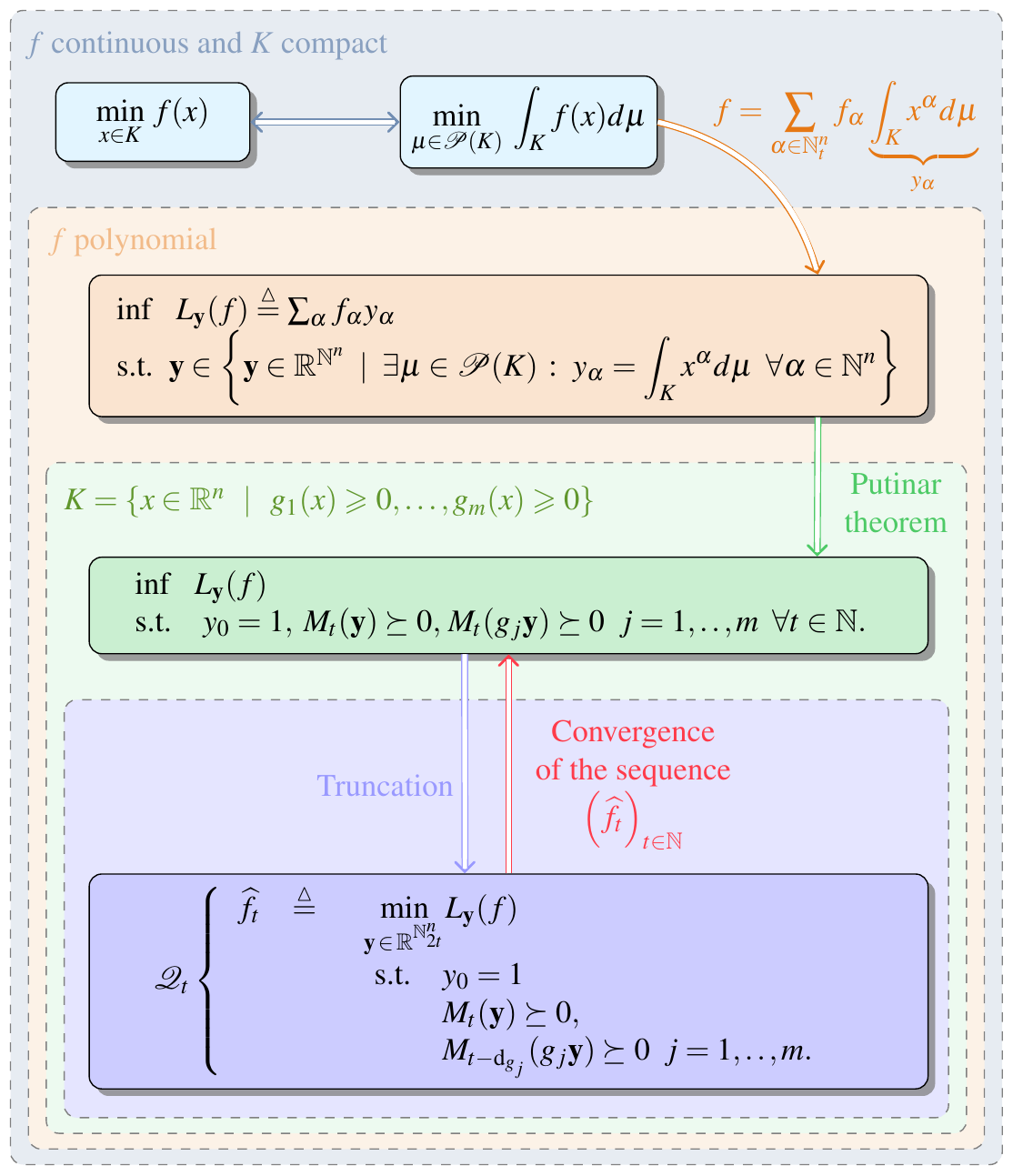}
\caption{Polynomial optimization process; see the main text for details.}
\end{figure*}

However, we do not know a priori at which relaxation order $t_0$ the convergence occurs. Practically, to detect whether the optimal value is attained, we can use conditions on the rank of the moment and localization matrices. Indeed, let $\y^\star \in \mathbb{R}^{\mathbb{N}_{2t}^n}$ be a solution of Problem~\eqref{def:fmomt} at a given relaxation order $t\geq \mathrm{d}_{\K}$, if:
\begin{eqnarray}\label{def:rang-conditions}
		\mathrm{rank}(M_{t}(\mathbf{y}^\star))=\mathrm{rank}(M_{t-\mathrm{d}_K}(\mathbf{y}^\star))
\end{eqnarray} 
then $\widehat{f}_t  = \widehat{f}$. In particular, if $\mathrm{rank}(M_{t}(\mathbf{y}^\star))=1\label{def:rang-conditions-1}$ then condition~\eqref{def:rang-conditions} is satisfied. Moreover, if these rank conditions are satisfied, then we can use numerical linear algebra to extract $\mathrm{rank}(M_{t}(\mathbf{y}^\star))$ global optima for Problem~\eqref{eq:prob_opti_globale}. We do not describe the algorithm in this article, but the reader can refer to~\cite[Sections 4.3]{LasserreBook:10} for more advanced information.
Figure~\ref{fig:gloptipoly} summarizes the optimization process.

A Matlab interface called GloptiPoly~\cite{gloptipoly03} has been designed to construct Lasserre's LMI relaxations in a format understandable by any SDP solver interfaced via YALMIP~\cite{YALMIP}. It can be used to construct an LMI relaxation~\eqref{def:fmomt} of a given order corresponding to a polynomial optimization problem~\eqref{eq:prob_opti_globale} with given polynomial data entered symbolically. A numerical algorithm is implemented in GloptiPoly to detect global optimality of an LMI relaxation, using the rank tests~\eqref{def:rang-conditions}. The algorithm also extracts numerically the global optima from the moment matrix.  
This approach has been successfully applied to globally solve various polynomial optimization problems (see \cite{LasserreBook:10} for an overview of results and applications). In computer vision this approach was first introduced in \cite{kahl:05} and used in \cite{KahlHenrion:07}.

\subsection{Application to Fundamental Matrix Estimation}
An important feature of our approach is that the singularity constraint can be directly satisfied by our solution and that we do {\it not} need an initial estimate
as in other methods. Nevertheless, the problem being homogeneous, an additional
normalization constraint is needed to avoid the trivial solution
$\Mat{F}=\textbf{0}$. This is generally done by setting one of the coefficients of the $\Mat{F}$ matrix to $\textbf{1}$.  However with such a normalization
there is no guarantee that the feasible set is compact which is 
a necessary condition for the convergence of our polynomial optimization method. Moreover, this normalisation excludes a priori some geometric configurations.
Thus, to guarantee compactness of the feasible set and to
avoid the trivial solution, we include the additional normalization constraint $\|\Mat{F}\|^2=1$.

Alternatively, the singularity constraint can be inferred by
parameterizing the $\Mat{F}$ matrix using one or two epipoles. For instance, in [10] the
authors have tried to use polynomial optimization to estimate the $\Mat{F}$
matrix, parametrized by using a single epipole.
However such an approach has several drawbacks.
Firstly, the coefficients of $\Mat{F}$ are not bounded and the convergence of
the method is not guaranteed.  Secondly, using the epipole explicitly
increases the degree of the polynomial criterion and consequently, the size of the
corresponding relaxations in the hierarchy. This results 
in a significant increase of the computational time.  Finally, the
choice of this parametrization is arbitrary and does not cover all camera configurations.

Our method is summarized in Algorithm~\ref{globalAlgorithm} below. Its main features are:
\begin{itemize}
\item In contrast with \cite{kahl:05,KahlHenrion:07}
the optimization problem is formulated 
with an {\it explicit} Frobenius norm constraint on the decision variables. This
enforces compactness of the feasibility set which is included
in the Euclidean ball of radius $1$. We have observed that enforcing
this Frobenius norm constraint has a dramatic influence on the overall numerical
behavior of the SDP solver, especially with respect to convergence and extraction of
global minimizers.
\item We have chosen the SDPT3 solver \cite{sdpt3,sdpt3UserGuide} since our experiments
revealed that for our problem it was the most efficient and
reliable solver.
\item We force the interior-point algorithm to increase the accuracy
as much as possible, overruling the default parameter set in SDPT3. Then the solver runs as long as it can make progress.
\item Generically, polynomial optimization problems on a compact set
have a unique global minimizer; this is confirmed in our numerical experiments
as the moment matrix has almost
always rank-one (which certifies global optimality) at the second SDP
relaxation of the hierarchy. In some (very few) cases, due to the numerical extraction, the global minimum is not fully accurate but yet largely satisfactory.
\end{itemize}

\begin{algorithm}
	\begin{algorithmic}[1]
		\REQUIRE Matched points $(\PtTwoD{q}_i,\PtTwoD{q}_i')$, $i=1,\dots,n$
		\STATE 	Create the cost function $Crit= \sum_{i=1}^n \left( \Transpose{\PtTwoD{q}_i'}\Mat{F}\PtTwoD{q}_i \right) ^2$:
		\begin{verbatim} mpol('F',3,3);
			 for k = 1:size(q1)
            			n(k) = (q2'*F*q1)^2;
			 end;
			 Crit = sum(n);
		\end{verbatim}
		\STATE Create the constraints $\det(\Mat{F})=0$ and $\lVert \Mat{F} \rVert^2 =1$:
			\begin{verbatim} K_det = det(F) == 0; K_fro = trace(F*F') == 1; \end{verbatim}
		\STATE Fix the accuracy of the solver to 0, then the solver runs as long as it can make progress:
		 \begin{verbatim} pars.eps = 0; mset(pars);\end{verbatim}
		\STATE Change the default SDP solver to SDPT3: \begin{verbatim}mset('yalmip',true); mset(sdpsettings('solver','sdpt3'));\end{verbatim}
		\STATE Form the second LMI relaxation of the problem:\begin{verbatim}P = msdp(min(crit),K_det,K_fro,2);\end{verbatim}
		\STATE Solve the second LMI relaxation:\begin{verbatim} msol(P);\end{verbatim}
	\end{algorithmic}
	\caption{\label{globalAlgorithm} Polynomial optimization for fundamental matrix estimation}
\end{algorithm}


\section{Experimental Results}
This section presents results obtained by the test procedure presented below with the 8-point method and our global method. First, criteria to evaluate the performance of a fundamental matrix estimate are described. Next, the evaluation methodology is detailed. Experiments were then carried out on synthetic data to test the sensitivity to noise and to the number of point matches. Finally, experiments on real data were performed to confirm previous results and to study the influence of the type ofmotion between the two images.

\subsection{Evaluation Criteria}
Various evaluation criteria were proposed in the literature~\cite{Zhang:98} to evaluate the quality of a fundamental matrix estimate. Driven by practice, a fundamental matrix estimate $\Mat{F}$ is evaluated with respect to the behavior of its subsequent refinement by projective bundle adjustment. 
Bundle Adjustment from two uncalibrated views is described as a minimization problem. The cost function is the RMS reprojection errors. The unknowns are the $3$D points $\Vct{Q}_i, \; i=1,\dots,n$ and the projection matrices $\Mat{P}$ and $\Mat{P}'$. The criteria we use are:
\begin{enumerate}
 \item The initial reprojection error written $e_{\text{Init}}(\Mat{F})$.
 \item The final reprojection error $e_{\text{BA}}(\Mat{F})$.
 \item The number of iterations taken by Bundle Adjustment to converge, $\text{Iter}(\Mat{F})$. 
\end{enumerate}
These three criteria assess whether the estimates provided by the two methods, denoted by $\Mat{F}_{8pt}$ and $\Mat{F}_{Gp}$, are in a `good' basin of attraction. Indeed, the number of iterations gives an indication of the distance between the estimate and the optimum while $e_{\text{BA}}(\Mat{F})$ gives an indication on the quality of the optimum.

\subsection{Evaluation Method}
The following algorithm summarizes our evaluation method:
\begin{enumerate}
  \item[] {\bf Evaluate}($\Mat{F}$)
  \item {\bf Inputs:} fundamental matrix estimate $\Mat{F}$, $n$ point matches $(\PtTwoD{q}_i,\PtTwoD{q}_i')$, $i=1,\dots,n$
  \item Form initial projective cameras~\cite{Luong:1996}:
  \item[]~~~~Find the second epipole from $\Transpose{\Mat{F}}\PtTwoD{e}'\sim \Vct{0}_{(3\times1)}$
  \item[]~~~~Find the canonical plane homography $\Mat{H}^*\sim [\PtTwoD{e}]_{\times}\Mat{F}$
  \item[]~~~~Set $\Mat{P} \sim \left[\Mat{I}_{(3\times3)} \; \Vct{0}_{(3\times1)}\right]$ and $\Mat{P}' \sim \left[\Mat{H}^* \; \PtTwoD{e}'\right]$
  \item Form initial 3D points~\cite{hartleysturmcviu}:
  \item[]~~~~Triangulate each point independently by minimizing the reprojection error
  \item Compute $e_{\text{Init}}(\Mat{F})$
  \item Run two-view uncalibrated bundle adjustment~\cite{Hartley-Zisserman:03}
  \item Compute $e_{\text{BA}}(\Mat{F})$ and $\text{Iter}(\Mat{F})$
  \item {\bf Outputs:} $e_{\text{Init}}(\Mat{F})$, $e_{\text{BA}}(\Mat{F})$ and $\text{Iter}(\Mat{F})$
\end{enumerate}

\subsection{Experiments on Simulated Data}
\subsubsection{Simulation Procedure}
\begin{figure*}[!htp]
  \centering
 \includegraphics[height = 1.5in]{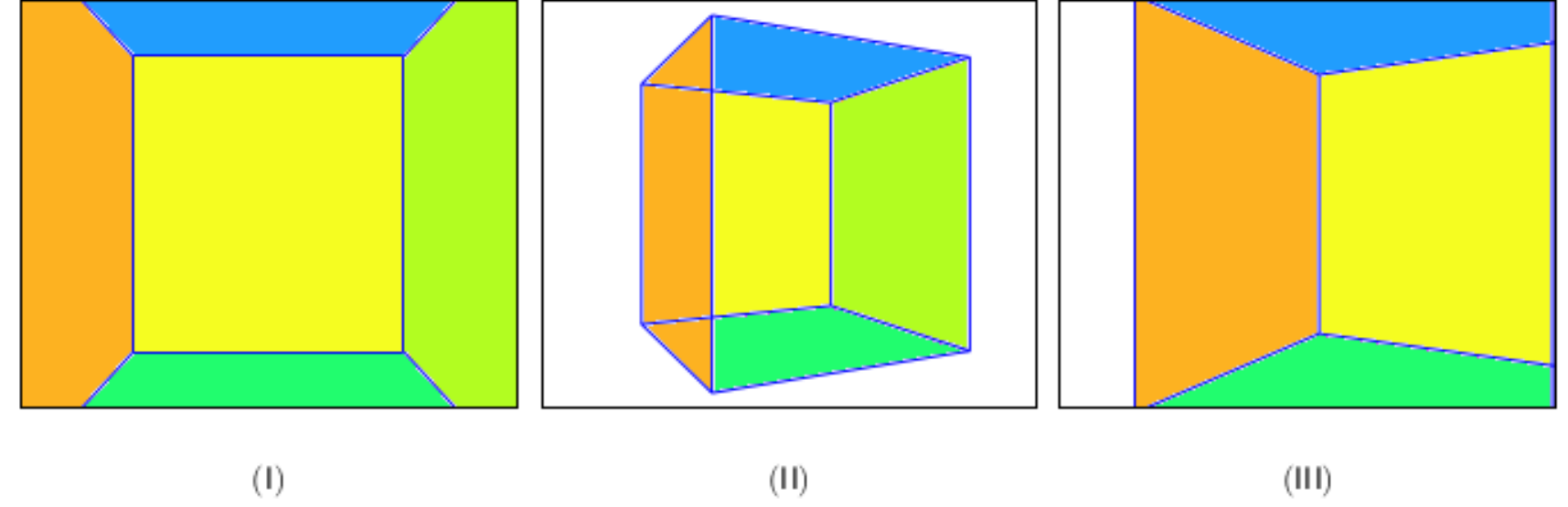}
\caption{\label{figure1}Projection of the cube in the camera on initial position $\mathsf{(I)}$ and in the camera after applying the rigid transformation $\left[\Mat{R}_1 \; t_1\right]$ $\mathsf{(II)}$ and the rigid transformation $\left[\Mat{R}_2\; t_2\right]$ $\mathsf{(III)}$. }
\end{figure*}

For each simulation series and for each parameter of interest (noise, number of points and number of motions), the
same methodology is applied with the following four steps:
\begin{enumerate}
\item For two given motions  between two successives images ($\left[\Mat{R}_k \; t_k \right]$) and for a given matrix $\Mat{K}$ of internal parameters, a set of 3D points $(\Vct{Q}_i)_{i}, \;i=1,\ldots,n$ is generated and two projection matrices $\Mat{P}$ and $\Mat{P}'$ are defined. In practice, the rotations matrices, $\Mat{R}_1$ and $\Mat{R}_2$, of two motions are defined by:
\begin{equation}\label{def:R1&R2}
\Mat{R}_k\stackrel{\vartriangle}{=}
	\left[
			\begin{array}{ccc}
					\cos(\theta_k)  &  0   & \sin(\theta_k)\\
					0               &  1   &  0\\
				       -\sin(\theta_k)  &  0   &  \cos(\theta_k) \\
			\end{array}
	\right]
\text{ with }
	\left\lbrace
			\begin{array}{l}
					\theta_1  =  \dfrac{\pi}{3} \\
					\text{ and }\\
				        \theta_2  = \dfrac{\pi}{6} \\
			\end{array}
	\right.
\end{equation}
and their translation vectors by $t_1=(20,0,5)^\top$ and $t_2=(6,0,0)^\top$. These matrices are chosen such that $\left[ \Mat{R}_1, t_1 \right] $ is a large movement and  $\left[ \Mat{R}_2, t_2 \right] $ is a small movement (see Figure~\ref{figure1}). We simulated points lying in a cube with 10 meter side length. The first camera looks at the center of the cube and it is located 15 meters from the center of the cube. The focal length of the camera is 700 pixels and the resolution is $640\times480$ pixels. 
\item Thanks to projection matrices $\Mat{P}=\Mat{K}\left[\Mat{R}_1 ,\Vct{t}_1\right]$ and $\Mat{P}'=\Mat{K}\left[\Mat{R}_2 ,\Vct{t}_2\right]$, the set of 3D points 
$(\PtThreeD{Q}_i)_{i}$ is projected into the two images as $(\PtTwoD{q}_i,\PtTwoD{q}'_i)_{i}$.
At each of their pixel coordinates, a centered Gaussian noise  with a variance $\sigma^2$ is added. In order to have statistical evidence, the results are averaged over 100 trials.
\item The resulting noisy points $(\widetilde{\PtTwoD{q}_i},\widetilde{\PtTwoD{q}'_i})_{i}$ 
are used to estimate $\Mat{F}$ by our method $\Mat{F}_{Gp}$ and the reference 8-point method $\Mat{F}_{8pt}$.
\item Finally, via our evaluation procedure we evaluate the estimation error with respect to the noise standard deviation $\sigma$ and the number of points $n$.
\end{enumerate}

\subsubsection{Sensitivity to Noise}
\begin{figure*}[!htp]
  \centering
   \includegraphics[height = 3.5in]{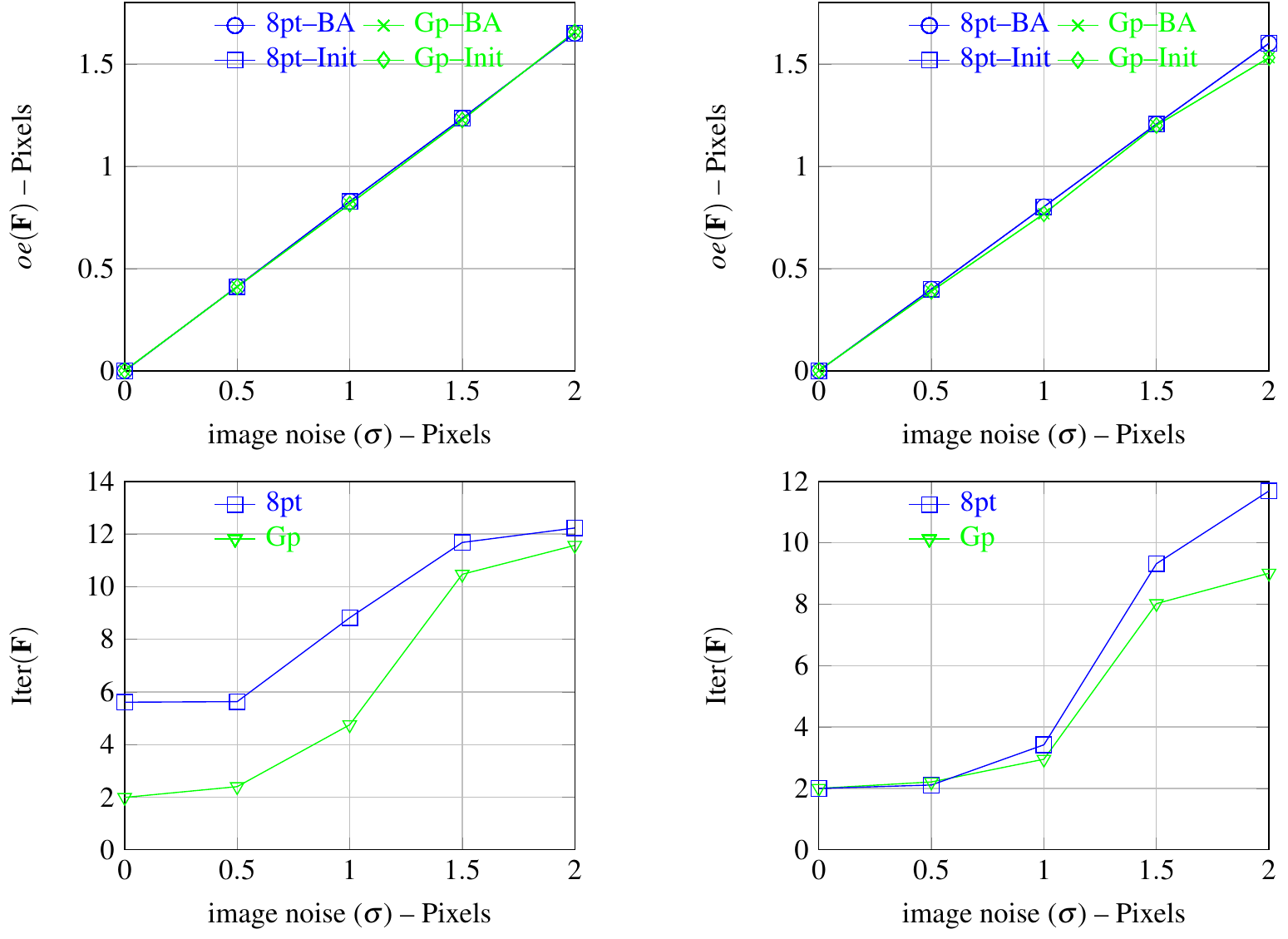}
\caption{\label{figure2} For two movements, $\left[\Mat{R}_1 \; t_1 \right]$ (left column) and $\left[\Mat{R}_2 \; t_2 \right]$ (right column), reprojection errors and number or iterations measured against image noise.}
\end{figure*}

We tested in two simulation series the influence of $\sigma$
 ranging from 0 to 2 pixels. The number of simulated points is $50$. The first (resp. second) simulation series is based on the first motion $\left[\Mat{R}_1 \; t_1 \right]$
(resp. the second motion $\left[\Mat{R}_2 \; t_2 \right]$). Figure \ref{figure2} gathers the influence of noise on the evaluation criteria. The first line shows the reproduction errors before, $e_{\text{Init}}(\Mat{F})$, and after $e_{\text{BA}}(\Mat{F})$ refinement through Bundle Adjustment with respect to the  noise standard deviation. 
The second line shows the number of iterations $\text{Iter}(\Mat{F})$ of the Bundle Adjustment versus the  noise standard deviation. The first (resp. second) row concerns the first (resp. second) motion. 

For the two motions, re-projection errors, $e_{\text{Init}}(\Mat{F})$ or $e_{\text{BA}}(\Mat{F})$, increase with the same slope when the noise level increases. Notice that for both movements, the Bundle-Adjustment step does not improve the results. Indeed, the noise gaussian noise is added to the projections $(\PtTwoD{q}_i,\PtTwoD{q}'_i)_{i}$. So this is noise which in practice would be produced by the extraction points process. Thus the solution produced by the resolution of the linear system is very close to the optimum and does not need to be refined. The initial solution provided by the triangulation step is then very close to a local minimum of the Bundle Adjustment problem. Moreover, the variation of the errors of initial re-projection before ($8pt-Init$ and $Gp-Init$) and after ($8pt-BA$ and $Gp-BA$) Bundle Adjustment versus the noise standard deviation is linear. However, the number of iterations needed for convergence is different in the two methods. The initial estimate of the triangulation computed 
from $\Mat{F}_{Gp}$ is closer to the local minimum than that obtained from  $\Mat{F}_{8pt}$. For the first motion (large displacement between camera 1 and 2), the number of iterations of the global method (in green) remains smaller than
for the 8-point method (in blue) even though their difference seems to decrease when the noise level is high $(\sigma>1$). 
For a significant displacement the quality of the estimate $\Mat{F}$ by the global method remains better even though
the difference in quality diminishes with the noise level. Conversely, for the second motion (small displacement between the camera 1 and 2) both methods are equivalent since the difference in quality is only significant for a high level of noise ($\sigma>1$). This is logical as the movement is less important. As a conclusion, the 8-point method provides a solution equivalent to that obtained with the global method when the displacement is not too important. For more significant movements the provided solution is not so close even though still in the same basin of attraction of a local minimum.

\subsubsection{Influence of the Number of Points}
In this experiment, we kept the noise level constant with a standard deviation $\sigma^2=0.5$ pixels. We tested the influence of the number of matches $(\PtTwoD{q}_i,\PtTwoD{q}'_i)_{i}$ on the quality of the resulting estimate of $\Mat{F}$. The number of points $N$ varied from 10 to 100. Two simulation series are also carried out with the two motions.

\begin{figure*}[!htp]
  \centering
   \includegraphics[height = 3.5in]{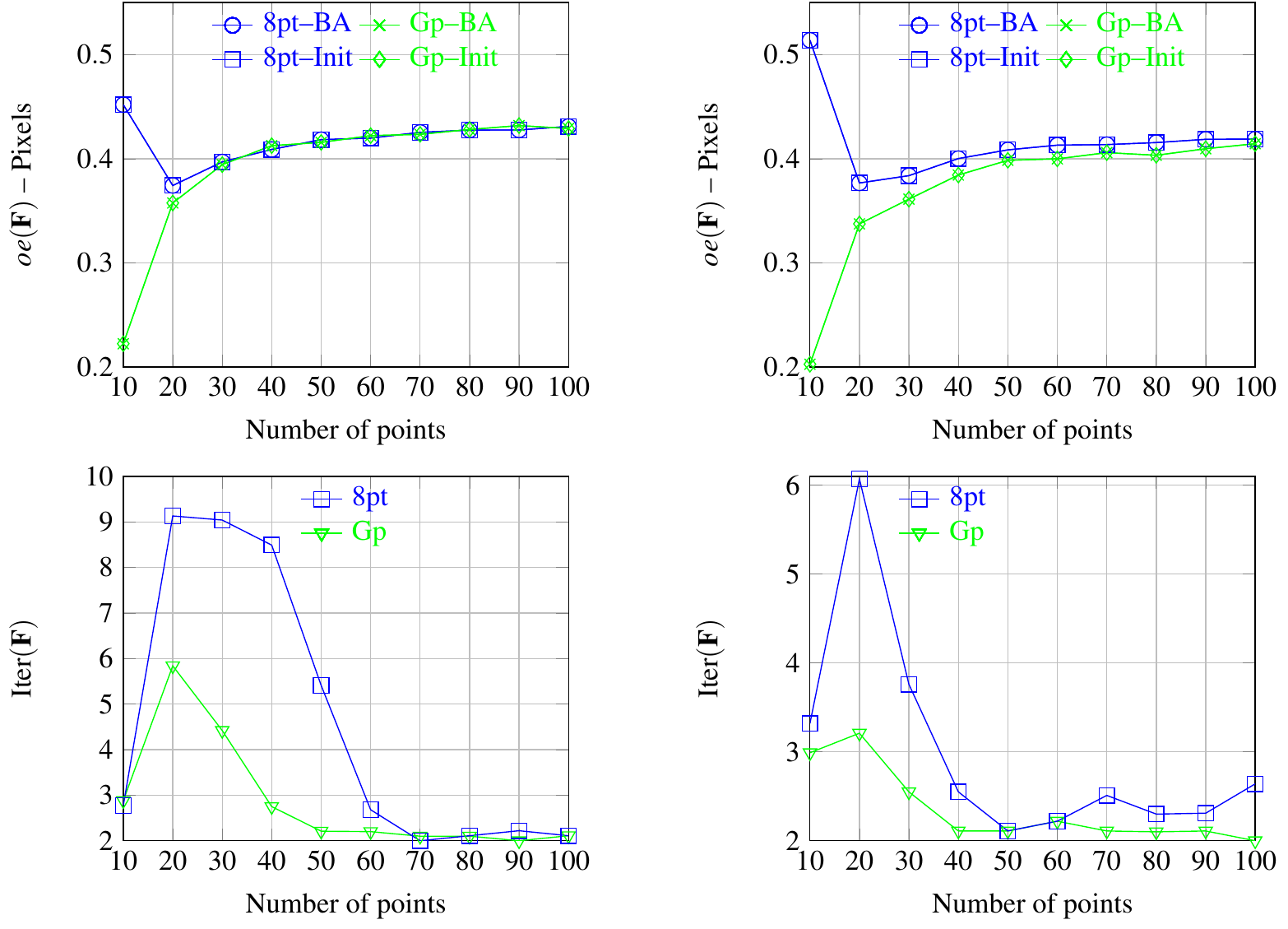}
\caption{\label{figure3}For two movements, $\left[\Mat{R}_1 \; t_1 \right]$ (left column) and $\left[\Mat{R}_2 \; t_2 \right]$ (right column), reprojection errors and number or iterations measured against number of points for a gaussian noise with a variance fixed to 0.5.}
\end{figure*}

Figure \ref{figure3} brings with the same organization the evaluation criteria. It displays the influence of the number of matches for estimating $\Mat{F}$ on the re-projection errors and on the number of iterations. For both motions and for a sufficiently high number of matches ($N>50$), re-projection errors, before and after refinement with Bundle Adjustment, or the number of iterations versus the  number of matches converge to the same asymptote. From a high number of matches, the initial estimate from triangulation computed with $\Mat{F}_{8pt}$ and with $\Mat{F}_{Gp}$ are both in the same basin of attraction for the Bundle Adjustment problem. However, for a number of matches smaller than $50$, the number of iterations to converge is smaller for given re-projection errors. The quality of the estimation by the global method seems better. The initial estimate from triangulation computed with $\Mat{F}_{8pt}$ goes away from the basin of convergence whereas the one computed with $\Mat{F}
_{Gp}$ remains in the basin. 

\subsubsection{Influence of the Number of Points with Wide Baseline}
\begin{figure*}[!htp]
  \centering
  \includegraphics[height = 1.78in]{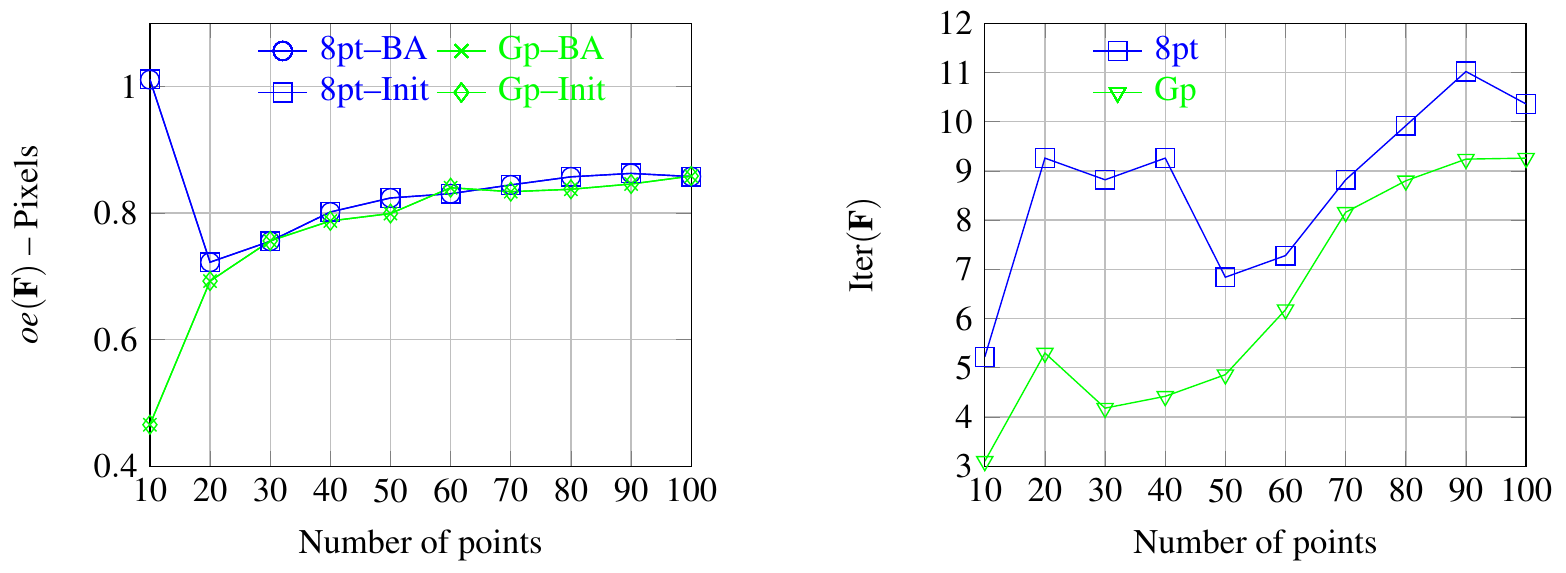}
\caption{\label{figure4} For the movement $\left[\Mat{R}_1 \; t_1 \right]$, reprojection errors and number or iterations measured against number of points for a gaussian noise with a variance fixed to 1 (left and right).}
\end{figure*}

In order to sustain the previous behavior, for a noise standard deviation of  $\sigma^2=1$ pixel and for the significant displacement $\left[\Mat{R}_1 \; t_1 \right]$, the influence of the number of matching points on the re-projection errors and on the number of iterations was tested. In this difficult context, Figure~\ref{figure4} demonstrates that the initial estimate computed from $\Mat{F}_{Gp}$ is always closer to the local minimum than that computed from $\Mat{F}_{8pt}$. No matter what is the number of matching points, the number of iterations needed to converge is always smaller.

As a conclusion, the quality of solutions obtained by both methods is almost identical when the movement is not too important, the number of matching points is sufficiently large, and the noise level is not too high. However, when one of these three parameters varies then the 8-point method lacks precision whereas the global method still allows Bundle Adjustment to convergence to the global minimum.
The 8-point method computes the projection of an unconstrained local minimizer on the feasible set whereas the global method provides a global minimizer of the constrained optimization problem. It is already surprising that even for good values of the three parameters the resulting solutions are not too far apart. But for worst values of the parameters it would be even more surprising.

\subsection{Experiments on Real Data}
The evaluation criteria remain the same, $e_{\text{Init}}(\Mat{F})$, $e_{\text{BA}}(\Mat{F})$ and $\text{Iter}(\Mat{F})$ 
and the computation time is added. Two experiments were carried out with two sets of images that illustrate different motions between two successive images.

\subsubsection{Experiment 1}
\begin{figure*}[!htp]
  \begin{center}
   \includegraphics[height = 4in]{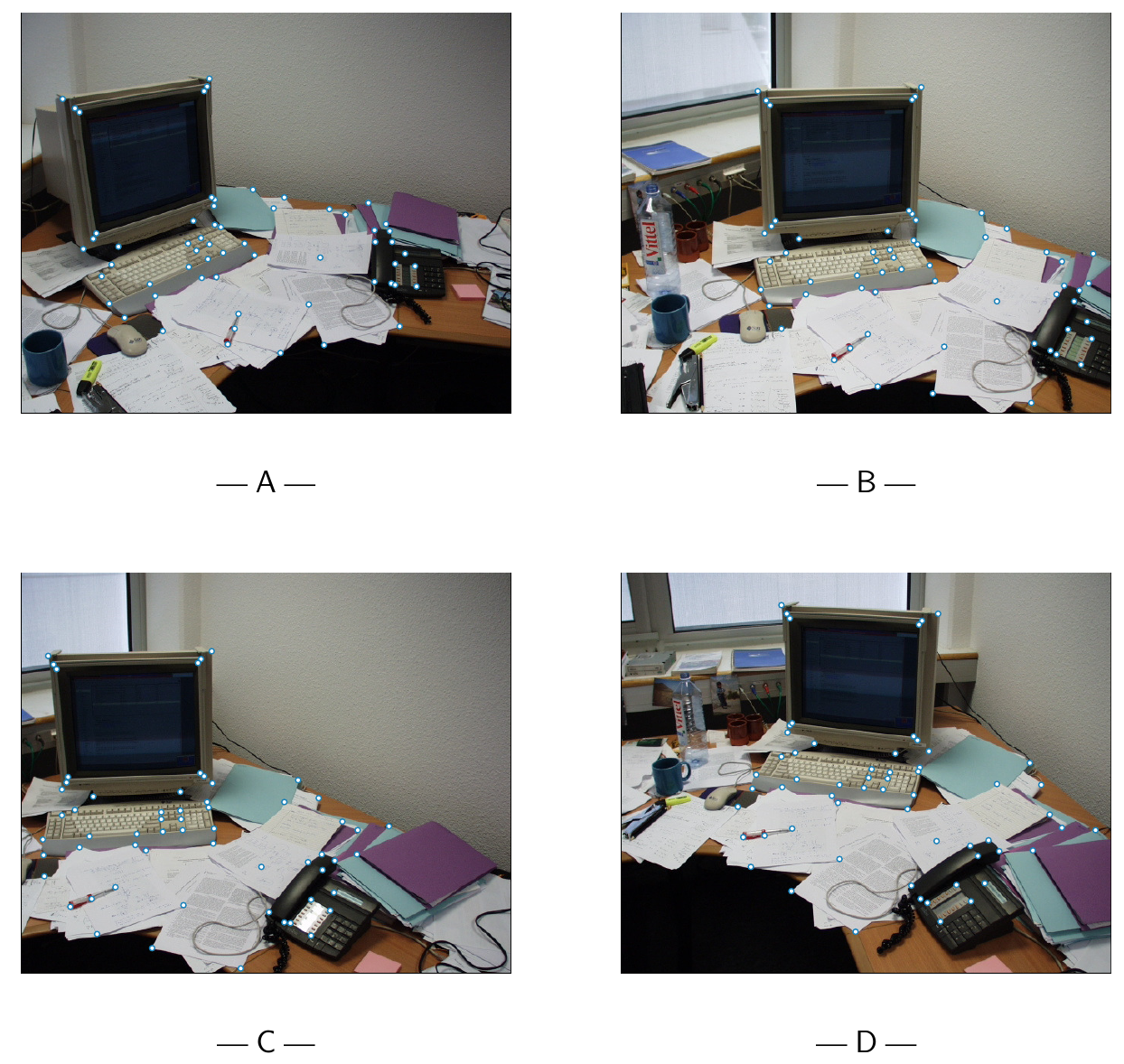}
  \end{center}
  \vspace*{0.25cm}
    \begin{center}
  \begin{tabular}{|c|c|c|c|c|c|c|c|c|c|c|c|c|}
	  \hline
	  \multicolumn{2}{|c|}{Views} &  \multicolumn{2}{c|}{Epipoles} &  \multicolumn{2}{c|}{$e_{\text{Init}}(\Mat{F})$} & 
	  \multicolumn{2}{c|}{$e_{\text{BA}}(\Mat{F})$} & \multicolumn{2}{c|}{$\text{Iter}(\Mat{F})$} & \multicolumn{2}{c|}{Time (s)}
	  \\
	  \cline{3-12}
	  \multicolumn{2}{|c|}{}  & $\PtTwoD{e}$ & $\PtTwoD{e}'$ &
	  $\Mat{F}_{8Pt}$ & $\Mat{F}_{Gp}$& 
	  $\Mat{F}_{8Pt}$ & $\Mat{F}_{Gp}$& 
	  $\Mat{F}_{8Pt}$ & $\Mat{F}_{Gp}$& 
	  $\Mat{F}_{8Pt}$ & $\Mat{F}_{Gp}$ 
	  \\
	  \hline
	  \hline
	  $\mathsf{A}$ & $\mathsf{B}$ & \sout{$\infty$} & $\infty$ & 0.597617  &0.65270 & 0.00252 & 0.00252 & 6 &	6 & 0.017	& 2.39 \\ 
	  \hline
	  $\mathsf{A}$ & $\mathsf{C}$ & \sout{$\infty$} & $\infty$ & 5.61506 & 5.61996 & 2.48258 & 0.00342 & 122 &	 175 & 0.017 & 1.98  \\
	  \hline
	  $\mathsf{A}$ & $\mathsf{D}$ & \sout{$\infty$} & \sout{$\infty$} & 21.0855 & 21.5848 &	4.74837 & 0.00344 &	105 &	30 & 0.017 & 2.12   \\
	  \hline
	  $\mathsf{B}$ & $\mathsf{C}$ & $\infty$  & $\infty$ & 2.49098 & 1.91136 & 0.00260 & 0.00260 & 	17 & 12 & 0.018 & 1.97  \\
	  \hline
	  $\mathsf{B}$ & $\mathsf{D}$ & $\infty$ &\sout{$\infty$} & 22.0071 & 23.6253 &	0.00268 & 0.00268 & 122 &	81 & 0.018 & 1.92 \\ 
	  \hline
	  $\mathsf{C}$ & $\mathsf{D}$ & $\infty$ &\sout{$\infty$} & 28.6586 & 28.6174 & 16.6507 & 0.25921 & 39 & 1001 & 0.018 & 2.1 \\
	  \hline
  \end{tabular}
\end{center}
\caption{\label{Table:AbDesk} Reprojection Error before ($e_{\text{Init}}(\Mat{F})$) and after Bundle Adjustment ($e_{\text{BA}}(\Mat{F})$), Number of Iterations ($\text{Iter}(\Mat{F})$), and CPU time to compute $\Mat{F}$ (Time), obtained when combining pairs of images to obtain epipoles close to the images or toward infinity}
\end{figure*}

The first set of four images (see Table \ref{Table:AbDesk}) shows all possible epipolar configurations
(right or left epipole at infinity $\dots$). With four images, six motions between a pair of images are possible:
$\mathsf{A}-\mathsf{B}$, $\mathsf{A}-\mathsf{C}$, $\mathsf{A}-\mathsf{D}$, $\mathsf{B}-\mathsf{C}$,
$\mathsf{B}-\mathsf{D}$ and $\mathsf{C}-\mathsf{D}$. For every pair of images,  $60$ matches are available to compute an estimate of $\Mat{F}$. The values ​​of the evaluation criteria are summarized in table~\ref{Table:AbDesk}. No matter what pair of images is used, the re-projection errors and the number of iterations are almost always better when  $\Mat{F}_{Gp}$ is used as initial guess. In addition, for three motions ($\mathsf{A}-\mathsf{C}$, $\mathsf{A}-\mathsf{D}$ and $\mathsf{C}-\mathsf{D}$), in contrast with the initial guess $\Mat{F}_{Gp}$, the initial guess from the 8-point method is not in a better basin of attraction. This may explain why the initial re-projection errors $e_{\text{Init}}(\Mat{F})$ are sometimes larger for $\Mat{F}_{Gp}$ as the initial guess may be in a good basin of attraction but with a larger re-projection error. For the four motions $\mathsf{A}-\mathsf{B}$, $\mathsf{B}-\mathsf{C}$, $\mathsf{B}-\mathsf{D}$ and  $\mathsf{B}-\mathsf{D}$, both initializations are in the same basin 
of attraction 
but the number of iterations demonstrates that the initial guess from the global method is always closer to the local minimizer. Finally, even though the computation time of the latter is significantly larger than for the 8-point method, it still remains compatible with a practical use.

\subsubsection{Experiment 2}
\begin{figure*}[htp]
  \begin{center}
   \includegraphics[height = 1.375in]{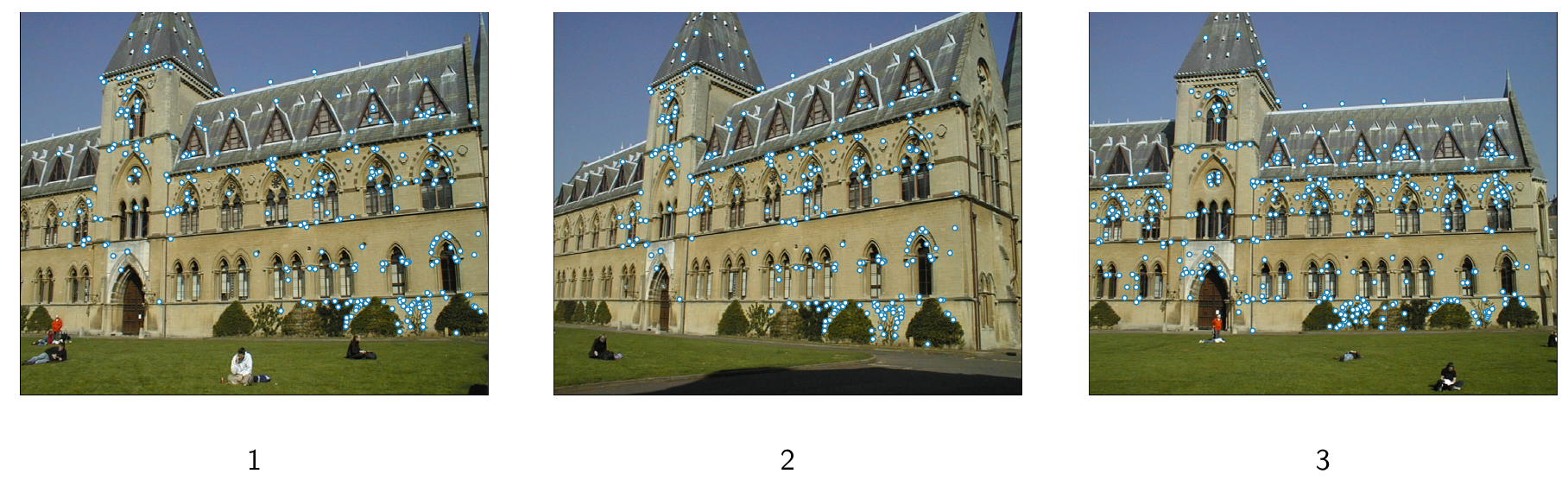}	\newline
   \\
  \begin{tabular}{|c|c|c|c|c|c|c|c|c|c|c|c|c|}
	  \hline
	  \multicolumn{2}{|c|}{Views} & \multicolumn{2}{c|}{nb Points} &   \multicolumn{2}{c|}{$e_{\text{Init}}(\Mat{F})$} & 
	  \multicolumn{2}{c|}{$e_{\text{BA}}(\Mat{F})$} & \multicolumn{2}{c|}{$\text{Iter}(\Mat{F})$} & \multicolumn{2}{c|}{Time (s)}
	  \\
	  \cline{5-12}
	  \multicolumn{2}{|c|}{} & \multicolumn{2}{c|}{}&
	  $\Mat{F}_{8Pt}$ & $\Mat{F}_{Gp}$& 
	  $\Mat{F}_{8Pt}$ & $\Mat{F}_{Gp}$& 
	  $\Mat{F}_{8Pt}$ & $\Mat{F}_{Gp}$& 
	  $\Mat{F}_{8Pt}$ & $\Mat{F}_{Gp}$ 
	  \\
	  \hline
	  \hline
	  $\mathsf{1}$ & $\mathsf{2}$ & \multicolumn{2}{c|}{310} & 18.27449 & 18.32095 & 11.39809 & 11.287591 & 65 & 57 & 0.035 & 2.05 \\ 
	  \hline
	  $\mathsf{2}$ & $\mathsf{3}$ & \multicolumn{2}{c|}{439} & 61.84817 & 63.79435 & 41.81854 & 40.2224 & 29 & 28 & 0.020 & 2.22  \\
	  \hline
	  $\mathsf{1}$ & $\mathsf{3}$ & \multicolumn{2}{c|}{439} & 42.67438 & 42.23867 & 28.65971 & 27.90128 & 40 & 32 & 0.043 & 2.10  \\
	  \hline
  \end{tabular}
\end{center}
\caption{\label{Table:Lib} Reprojection Error before ($e_{\text{Init}}(\Mat{F})$) and after Bundle Adjustment ($e_{\text{BA}}(\Mat{F})$), Number of Iterations ($\text{Iter}(\Mat{F})$), and CPU time to compute $\Mat{F}$ (Time), obtained when combining pairs of images of the \textit{Library} series}
\end{figure*}

\begin{figure*}[htp]
  \begin{center}
   \includegraphics[height = 1.375in]{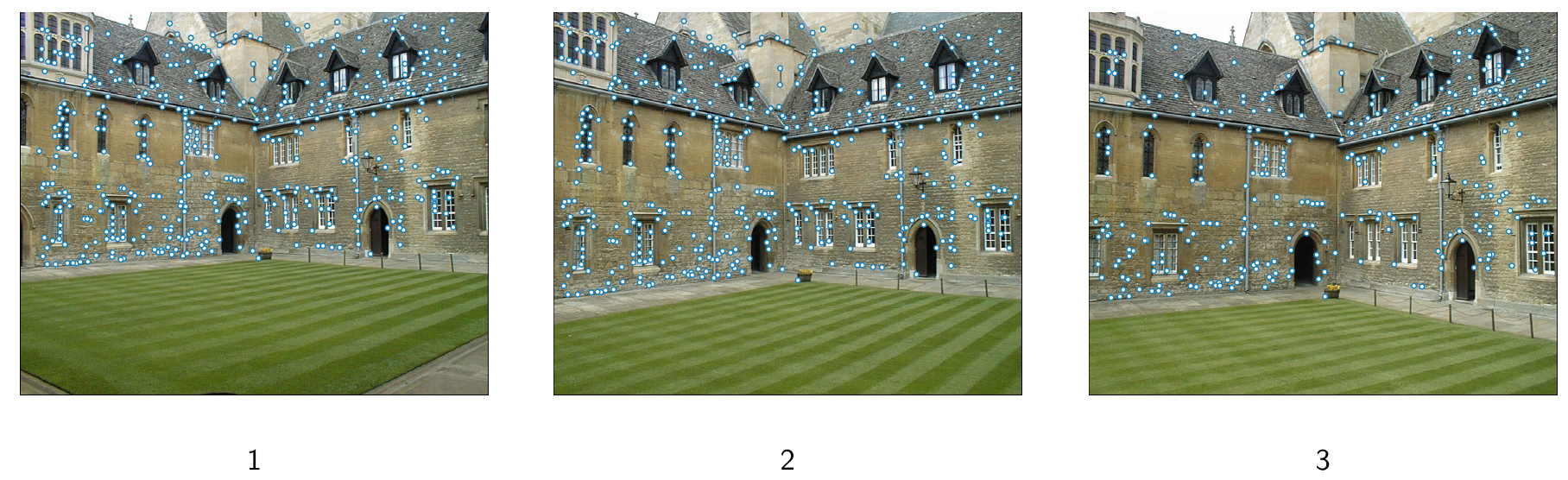}	\newline
   \\
  \begin{tabular}{|c|c|c|c|c|c|c|c|c|c|c|c|c|}
	  \hline
	  \multicolumn{2}{|c|}{Views} & \multicolumn{2}{c|}{nb Points} &   \multicolumn{2}{c|}{$e_{\text{Init}}(\Mat{F})$} & 
	  \multicolumn{2}{c|}{$e_{\text{BA}}(\Mat{F})$} & \multicolumn{2}{c|}{$\text{Iter}(\Mat{F})$} & \multicolumn{2}{c|}{Time (s)}
	  \\
	  \cline{5-12}
	  \multicolumn{2}{|c|}{} & \multicolumn{2}{c|}{}&
	  $\Mat{F}_{8Pt}$ & $\Mat{F}_{Gp}$& 
	  $\Mat{F}_{8Pt}$ & $\Mat{F}_{Gp}$& 
	  $\Mat{F}_{8Pt}$ & $\Mat{F}_{Gp}$& 
	  $\Mat{F}_{8Pt}$ & $\Mat{F}_{Gp}$ 
	  \\
	  \hline
	  \hline
	  $\mathsf{1}$ & $\mathsf{2}$ & \multicolumn{2}{c|}{485} & 15.1771 & 15.3126 & 9.8859 & 9.7360 & 77 & 55 & 0.048 & 2.02 \\ 
	  \hline
	  $\mathsf{2}$ & $\mathsf{3}$ & \multicolumn{2}{c|}{439} & 61.84817 & 63.79435 & 41.81854 & 40.2224 & 29 & 28 & 0.020 & 2.22  \\
	  \hline
	  $\mathsf{1}$ & $\mathsf{3}$ & \multicolumn{2}{c|}{384} & 51.3128 & 2.3690 & 7.5245 & 0.56180 & 5 & 52 & 0.040 & 2.71  \\
	  \hline
  \end{tabular}
\end{center}
\caption{\label{Table:Merton1} Reprojection Error before ($e_{\text{Init}}(\Mat{F})$) and after Bundle Adjustment ($e_{\text{BA}}(\Mat{F})$), Number of Iterations ($\text{Iter}(\Mat{F})$), and CPU time to compute $\Mat{F}$ (Time), obtained when combining pairs of images of the \textit{Merton1} series}
\end{figure*}

\begin{figure*}[htp]
  \begin{center}
   \includegraphics[height = 1.375in]{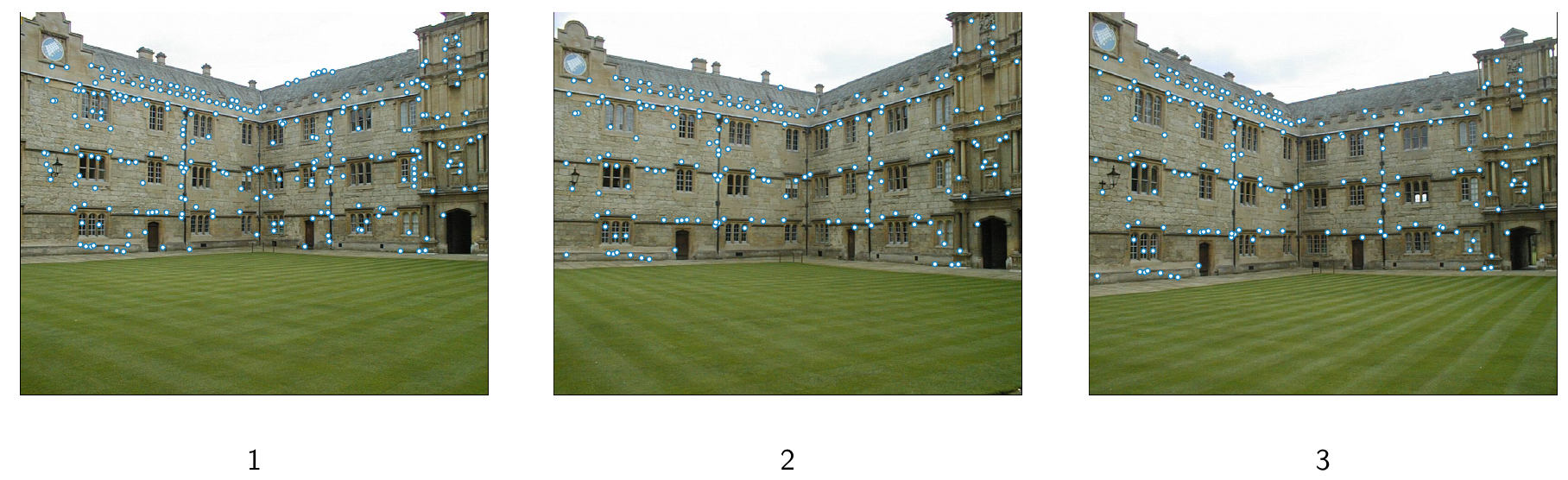}	\newline
   \\
  \begin{tabular}{|c|c|c|c|c|c|c|c|c|c|c|c|c|}
	  \hline
	  \multicolumn{2}{|c|}{Views} & \multicolumn{2}{c|}{nb Points} &   \multicolumn{2}{c|}{$e_{\text{Init}}(\Mat{F})$} & 
	  \multicolumn{2}{c|}{$e_{\text{BA}}(\Mat{F})$} & \multicolumn{2}{c|}{$\text{Iter}(\Mat{F})$} & \multicolumn{2}{c|}{Time (s)}
	  \\
	  \cline{5-12}
	  \multicolumn{2}{|c|}{} & \multicolumn{2}{c|}{}&
	  $\Mat{F}_{8Pt}$ & $\Mat{F}_{Gp}$& 
	  $\Mat{F}_{8Pt}$ & $\Mat{F}_{Gp}$& 
	  $\Mat{F}_{8Pt}$ & $\Mat{F}_{Gp}$& 
	  $\Mat{F}_{8Pt}$ & $\Mat{F}_{Gp}$ 
	  \\
	  \hline
	  \hline
	  $\mathsf{1}$ & $\mathsf{2}$ & \multicolumn{2}{c|}{345} & 24.7158 & 29.6530 & 14.0481 & 2.7599 & 69 & 42 & 0.037 & 2.91 \\ 
	  \hline
	  $\mathsf{2}$ & $\mathsf{3}$ & \multicolumn{2}{c|}{197} & 143.431 & 154.880 & 73.8776 & 72.8977 & 14 & 19 & 0.026 & 2.52  \\
	  \hline
	  $\mathsf{1}$ & $\mathsf{3}$ & \multicolumn{2}{c|}{270} & 59.8109 & 77.9734 & 30.3824 & 15.6367 & 38 & 32 & 0.031 & 2.37  \\
	  \hline
  \end{tabular}
\end{center}
\caption{\label{Table:Merton2} Reprojection Error before ($e_{\text{Init}}(\Mat{F})$) and after Bundle Adjustment ($e_{\text{BA}}(\Mat{F})$), Number of Iterations ($\text{Iter}(\Mat{F})$), and CPU time to compute $\Mat{F}$ (Time), obtained when combining pairs of images of the \textit{Merton2} series}
\end{figure*}

\begin{figure*}[htp]
  \begin{center}
   \includegraphics[height = 1.375in]{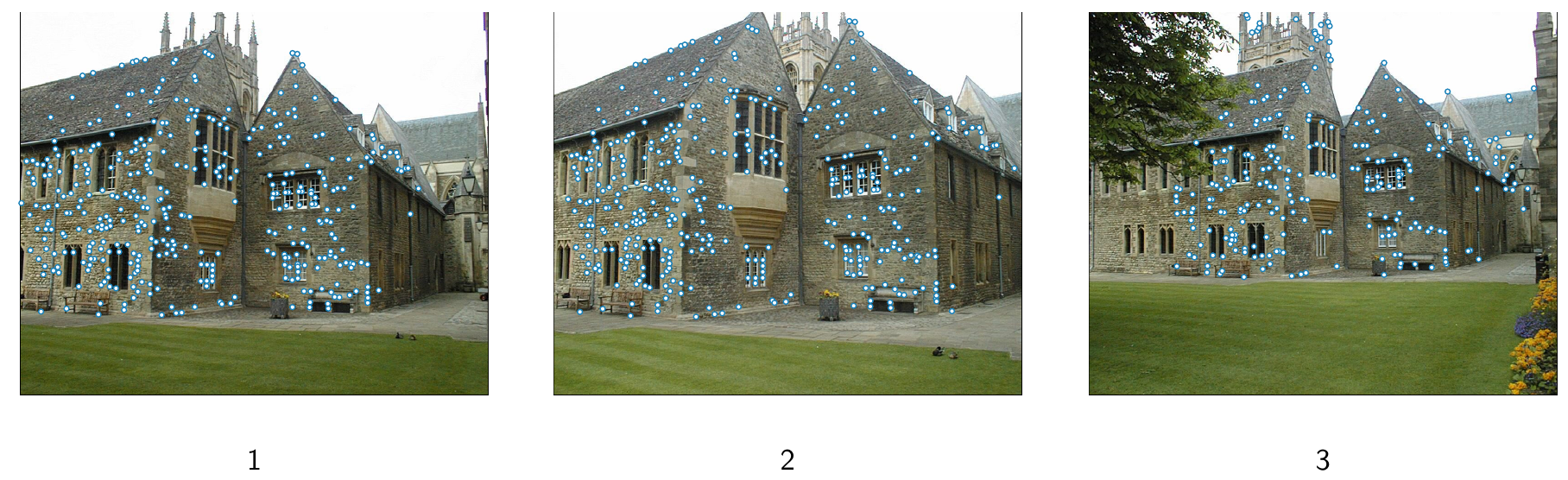}	\newline
   \\
  \begin{tabular}{|c|c|c|c|c|c|c|c|c|c|c|c|c|}
	  \hline
	  \multicolumn{2}{|c|}{Views} & \multicolumn{2}{c|}{nb Points} &   \multicolumn{2}{c|}{$e_{\text{Init}}(\Mat{F})$} & 
	  \multicolumn{2}{c|}{$e_{\text{BA}(\Mat{F})}$} & \multicolumn{2}{c|}{$\text{Iter}(\Mat{F})$} & \multicolumn{2}{c|}{Time (s)}
	  \\
	  \cline{5-12}
	  \multicolumn{2}{|c|}{} & \multicolumn{2}{c|}{}&
	  $\Mat{F}_{8Pt}$ & $\Mat{F}_{Gp}$& 
	  $\Mat{F}_{8Pt}$ & $\Mat{F}_{Gp}$& 
	  $\Mat{F}_{8Pt}$ & $\Mat{F}_{Gp}$& 
	  $\Mat{F}_{8Pt}$ & $\Mat{F}_{Gp}$ 
	  \\
	  \hline
	  \hline
	  $\mathsf{1}$ & $\mathsf{2}$ & \multicolumn{2}{c|}{400} & 62.7290 & 68.6030 & 52.6557 & 23.1740 & 5 & 32 & 0.041 & 3.03 \\ 
	  \hline
	  $\mathsf{2}$ & $\mathsf{3}$ & \multicolumn{2}{c|}{197} & 135.950 & 140.801 & 83.3365 & 76.8614 & 13 & 14 & 0.025 & 2.52  \\
	  \hline
	  $\mathsf{1}$ & $\mathsf{3}$ & \multicolumn{2}{c|}{264} & 116.7659 & 118.1754 & 24.7184 & 13.5496 & 9 & 11 & 0.032 & 2.39  \\
	  \hline
  \end{tabular}
\end{center}
\caption{\label{Table:Merton3} Reprojection Error before ($e_{\text{Init}}(\Mat{F})$) and after Bundle Adjustment ($e_{\text{BA}}(\Mat{F})$), Number of Iterations ($\text{Iter}(\Mat{F})$), and CPU time to compute $\Mat{F}$ (Time), obtained when combining pairs of images of the \textit{Merton3} series }
\end{figure*}

The second experiment compares the two methods on large motions. It is based on many series of images.  First we test our algorithm with the classic series \emph{Library}, \emph{Merton}, \emph{dinosaur} and \emph{house} that are available at \url{www.robots.ox.ac.uk/~vgg/data/data-mview.html}. For the set of three images of \emph{Library} and \emph{Merton} serie, Table~\ref{Table:Lib},~\ref{Table:Merton1},~\ref{Table:Merton2} and Table~\ref{Table:Merton3} demonstrate that the quality of the solution achieved by the global method is always better than with the 8-point method (in some cases both solutions are very close). 

\begin{figure*}[htp]
  \begin{center}
   \includegraphics[height = 1.375in]{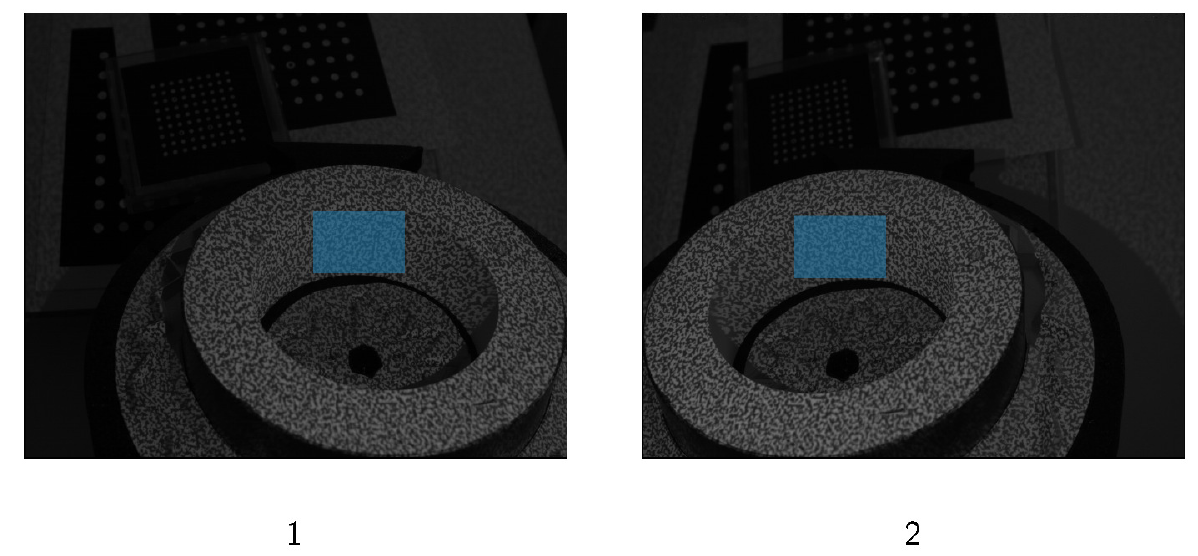}	\newline
   \\
  \begin{tabular}{|c|c|c|c|c|c|c|c|c|c|c|c|c|}
	  \hline
	  \multicolumn{2}{|c|}{Views} & \multicolumn{2}{c|}{nb Points} &   \multicolumn{2}{c|}{$e_{\text{Init}}(\Mat{F})$} & 
	  \multicolumn{2}{c|}{$e_{\text{BA}}(\Mat{F})$} & \multicolumn{2}{c|}{$\text{Iter}(\Mat{F})$} & \multicolumn{2}{c|}{Time (s)}
	  \\
	  \cline{5-12}
	  \multicolumn{2}{|c|}{} & \multicolumn{2}{c|}{}&
	  $\Mat{F}_{8Pt}$ & $\Mat{F}_{Gp}$& 
	  $\Mat{F}_{8Pt}$ & $\Mat{F}_{Gp}$& 
	  $\Mat{F}_{8Pt}$ & $\Mat{F}_{Gp}$& 
	  $\Mat{F}_{8Pt}$ & $\Mat{F}_{Gp}$ 
	  \\
	  \hline
	  \hline
	  $\mathsf{1}$ & $\mathsf{2}$ & \multicolumn{2}{c|}{6609} & 6.3142 & 6.4021 & 1.4701 & 0.2531 & 401 & 237 & 1.23 & 2.33 \\ 
	  \hline
  \end{tabular}
\end{center}
\caption{\label{Table:Cylinder} Reprojection Error before ($e_{\text{Init}}(\Mat{F})$) and after Bundle Adjustment ($e_{\text{BA}}(\Mat{F})$), Number of Iterations ($\text{Iter}(\Mat{F})$), and CPU time to compute $\Mat{F}$ (Time), obtained when combining pairs of images of the \textit{Cylinder} series. The matched points are located in blue bounding boxes.}
\end{figure*}

\begin{figure*}[htp]
  \begin{center}
   \includegraphics[height = 1.375in]{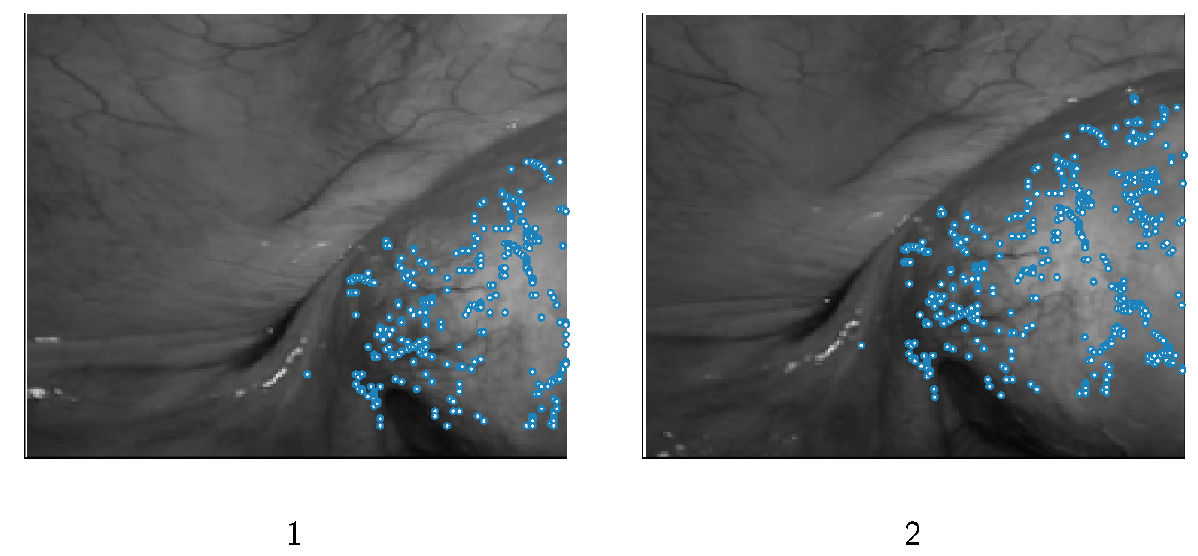}	\newline
   \\
  \begin{tabular}{|c|c|c|c|c|c|c|c|c|c|c|c|c|}
	  \hline
	  \multicolumn{2}{|c|}{Views} & \multicolumn{2}{c|}{nb Points} &   \multicolumn{2}{c|}{$e_{\text{Init}}(\Mat{F})$} & 
	  \multicolumn{2}{c|}{$e_{\text{BA}}(\Mat{F})$} & \multicolumn{2}{c|}{$\text{Iter}(\Mat{F})$} & \multicolumn{2}{c|}{Time (s)}
	  \\
	  \cline{5-12}
	  \multicolumn{2}{|c|}{} & \multicolumn{2}{c|}{}&
	  $\Mat{F}_{8Pt}$ & $\Mat{F}_{Gp}$& 
	  $\Mat{F}_{8Pt}$ & $\Mat{F}_{Gp}$& 
	  $\Mat{F}_{8Pt}$ & $\Mat{F}_{Gp}$& 
	  $\Mat{F}_{8Pt}$ & $\Mat{F}_{Gp}$ 
	  \\
	  \hline
	  \hline
	  $\mathsf{1}$ & $\mathsf{2}$ & \multicolumn{2}{c|}{730} & 5.6624 & 5.6366 & 3.8566 & 0.9386 & 163 & 401 & 0.07 & 1.69 \\ 
	  \hline
  \end{tabular}
\end{center}
\caption{\label{Table:Endoscope} Reprojection Error before ($e_{\text{Init}}(\Mat{F})$) and after Bundle Adjustment ($e_{\text{BA}}(\Mat{F})$), Number of Iterations ($\text{Iter}(\Mat{F})$), and CPU time to compute $\Mat{F}$ (Time), obtained when combining pairs of images of the \textit{Endoscope} series}
\end{figure*}

We also conducted the same tests on other pairs of images. For the first pair, we used images from a standard cylinder graciously provided by the company NOOMEO (\url{http://www.noomeo.eu/}). This cylinder is use to evaluate the accuracy of 3D reconstructions. Matched points are calculated with digital image correlation method. They are located in a window inside the cylinder. Thus, we have 6609 pairs $(\PtTwoD{q}_i,\PtTwoD{q}_i)_i$ matched to sub-pixel precision. Results are presented in the Table~\ref{Table:Cylinder}. We observe that the computation time of the 8-point method exceeds one second. This is due to the large number of matched points which leads to the resolution of a large linear system. However, as the points are precisely matched, this system is well conditioned. But the quality of the fundamental matrix estimated with the 8-point method is not sufficient to properly initialize the Bundle-Adjustment because the final re-projection error is 1.47 pixels. At the same time, even if the number of 
iterations is larger, our global method supplies a good estimation because the final re-projection error is 0.25 pixels. Furthermore, the calculation time remains constant in approximately 2 seconds. For the second pair, we use images taken by an endoscope. Table~\ref{Table:Endoscope} shows the results obtained on this difficult case. As for the previous example, we observe that the fundamental matrix estimated by our global method is good quality because the final error is 0.93 pixels. At the same time, Bundle Adjustement puts more iterations to converge on a less precise solution when we use $\Mat{F}_{8pt}$ to initialize it.

\begin{figure*}[htp]
\centering
      \includegraphics[height = 2.5in]{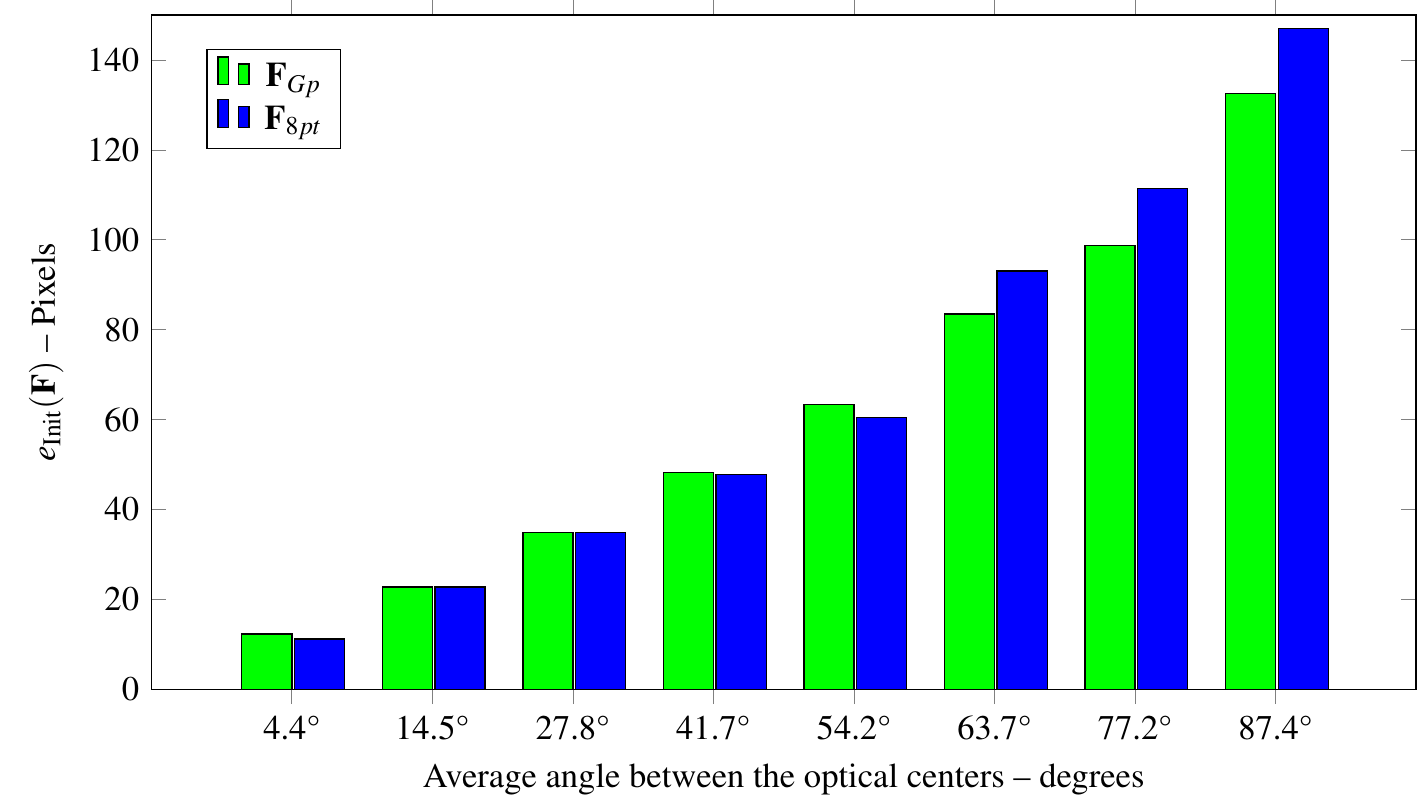}
\caption{\label{figure12} Initial re-projections errors measured against movement amplitude for the \textit{dinosaur} series}
\end{figure*}

\begin{figure*}[htp]
\centering
   \includegraphics[height = 2.5in]{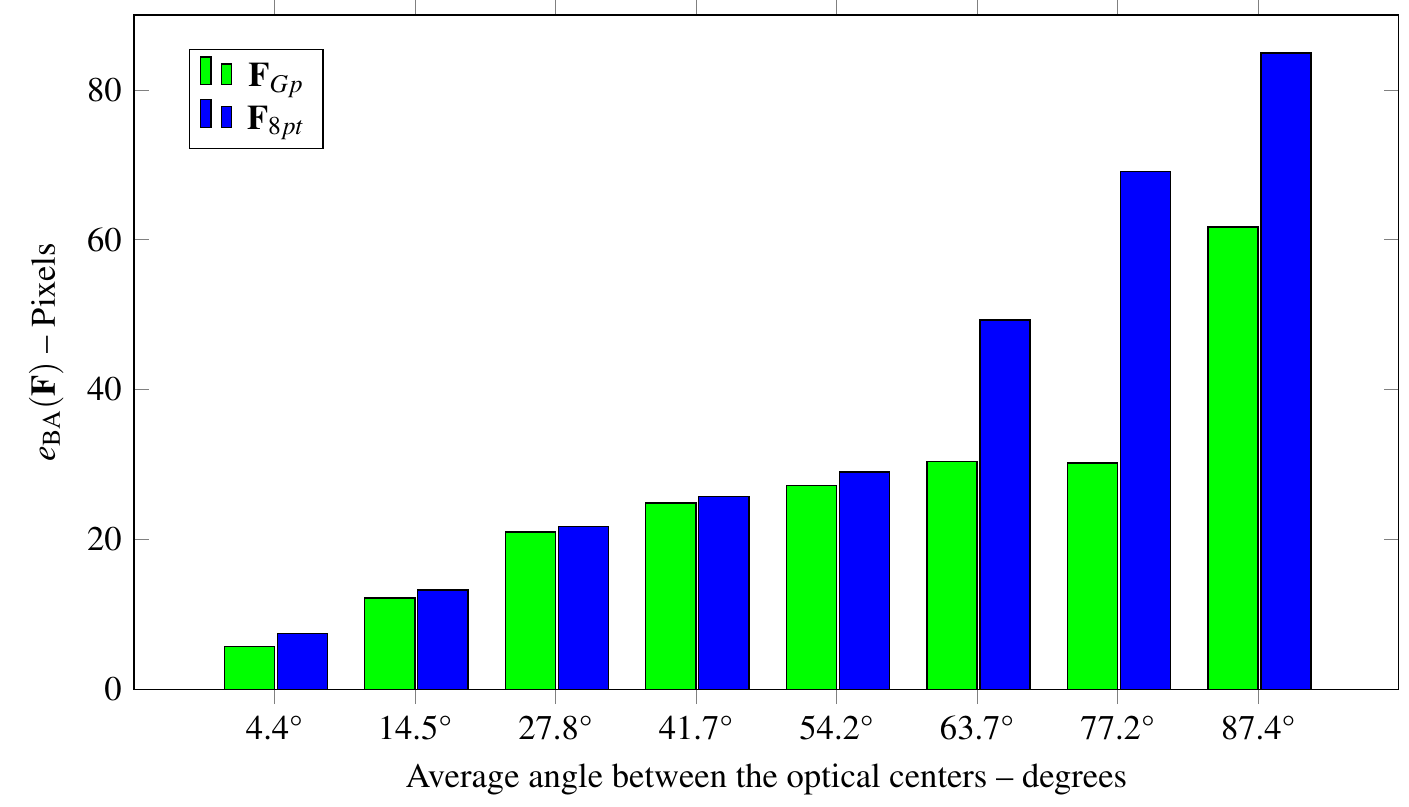}
\caption{\label{figure13} Final re-projections errors measured against movement amplitude for the \textit{dinosaur} series}
\end{figure*}

\begin{figure*}[htp]
\centering
   \includegraphics[height = 2.5in]{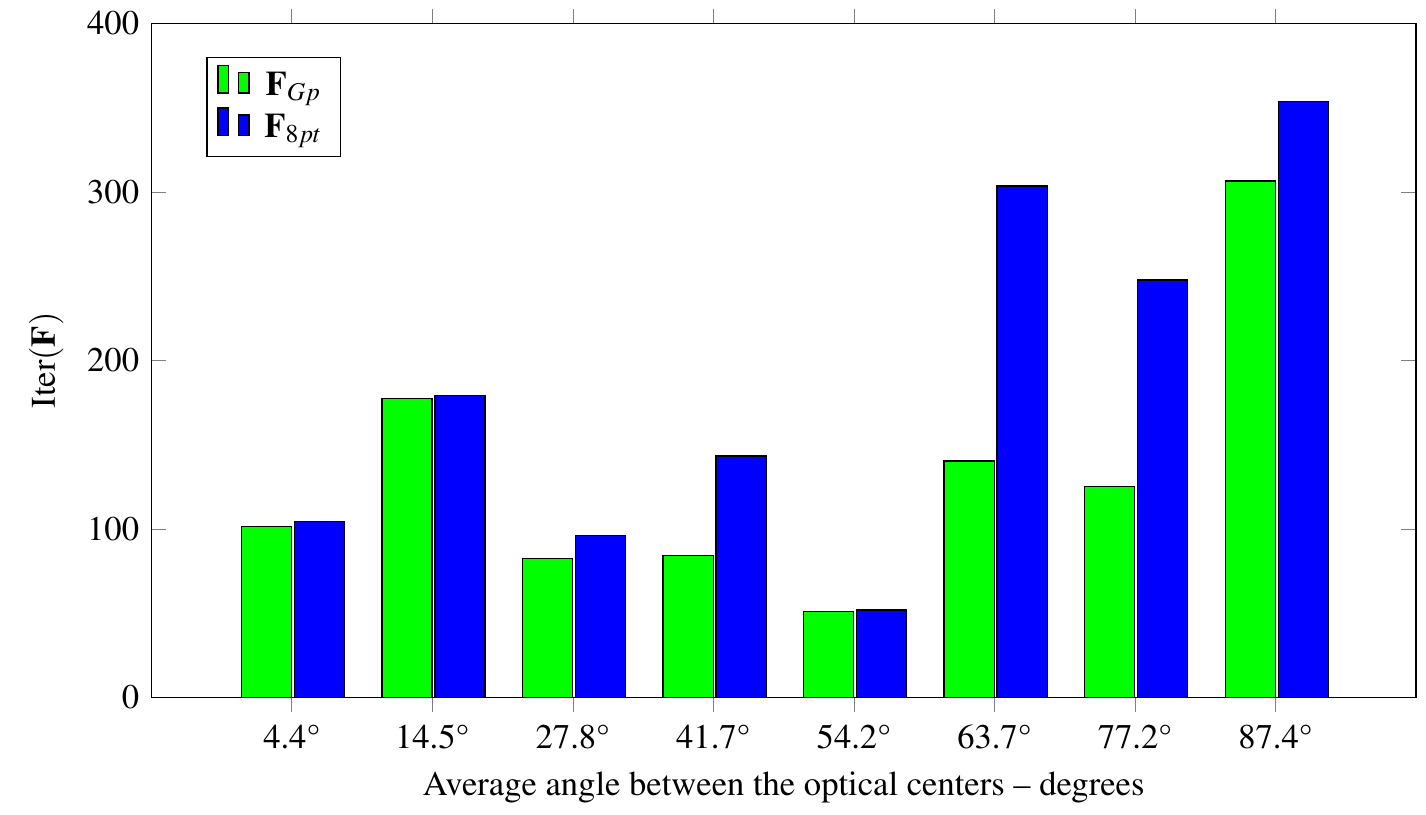}
\caption{\label{figure14} Number of iterations performed by Bundle-Adjustment to converge measured against movement amplitude for the \textit{dinosaur} series}
\end{figure*}

\begin{figure*}[htp]
\centering
   \includegraphics[height = 2.5in]{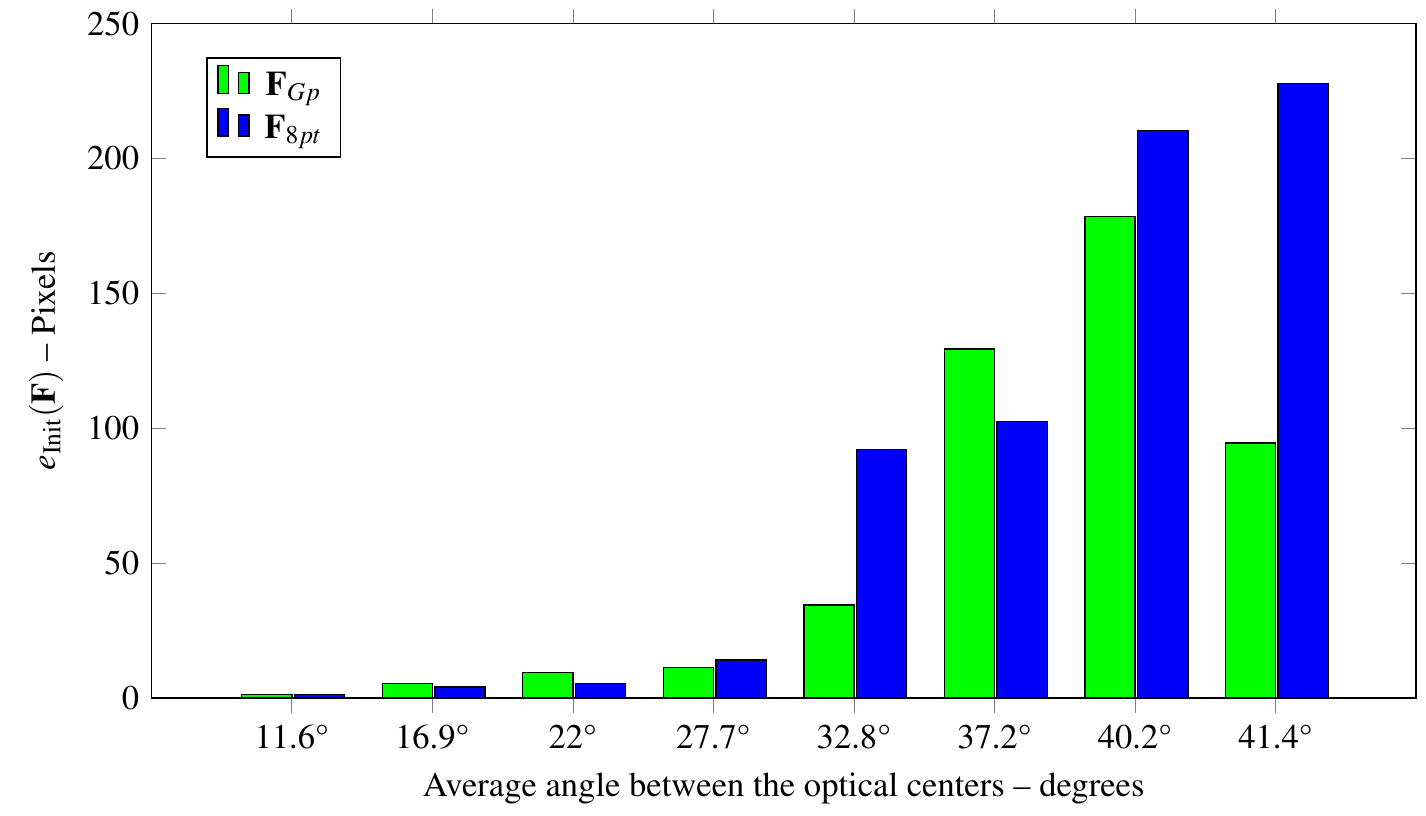}
\caption{\label{figure15} Initial re-projections errors measured against movement amplitude for the \textit{House} series }
\end{figure*}

\begin{figure*}[htp]
\centering
   \includegraphics[height = 2.5in]{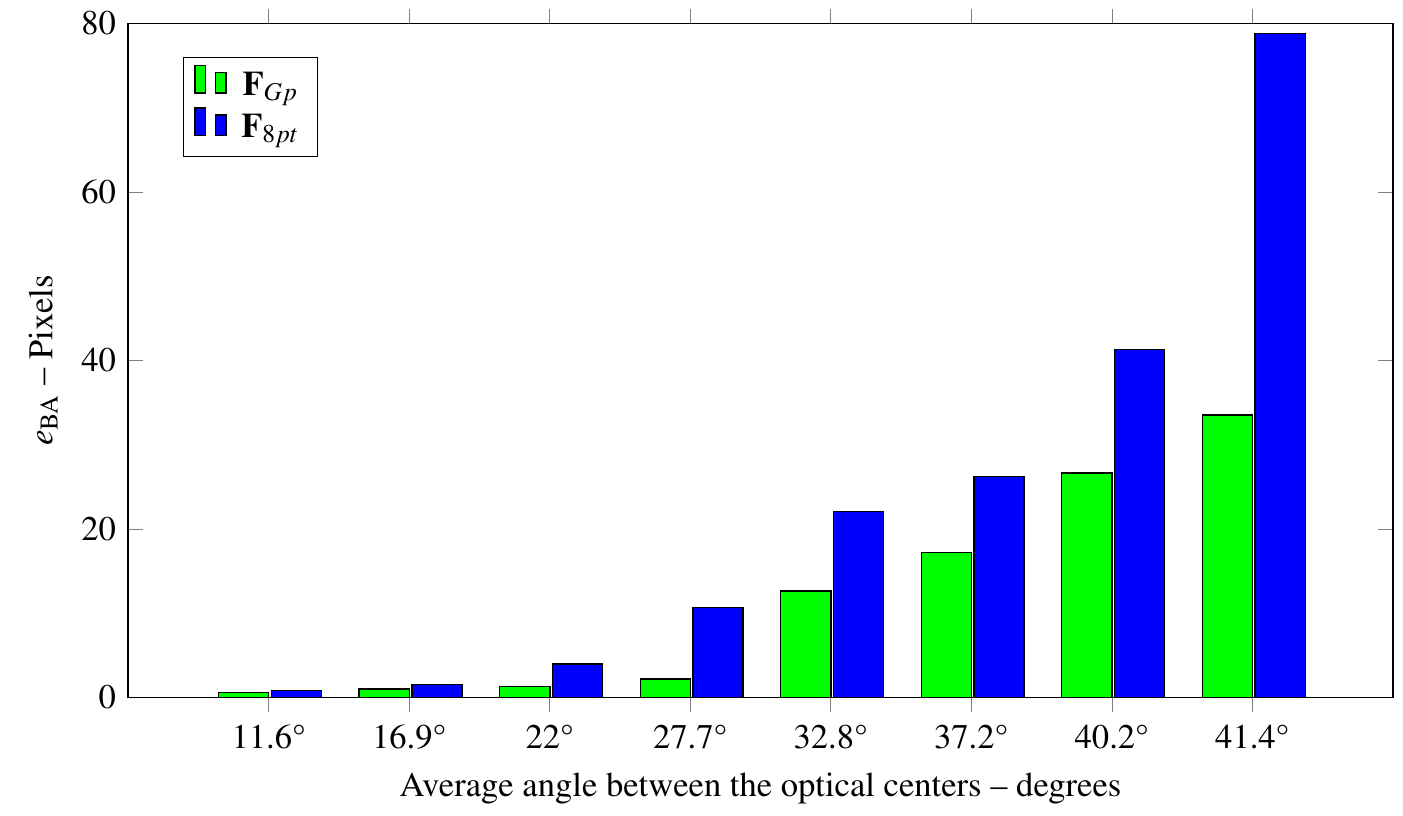}
\caption{\label{figure16} Final re-projections errors measured against movement amplitude for the \textit{House} series }
\end{figure*}

\begin{figure*}[htp]
\centering
   \includegraphics[height = 2.5in]{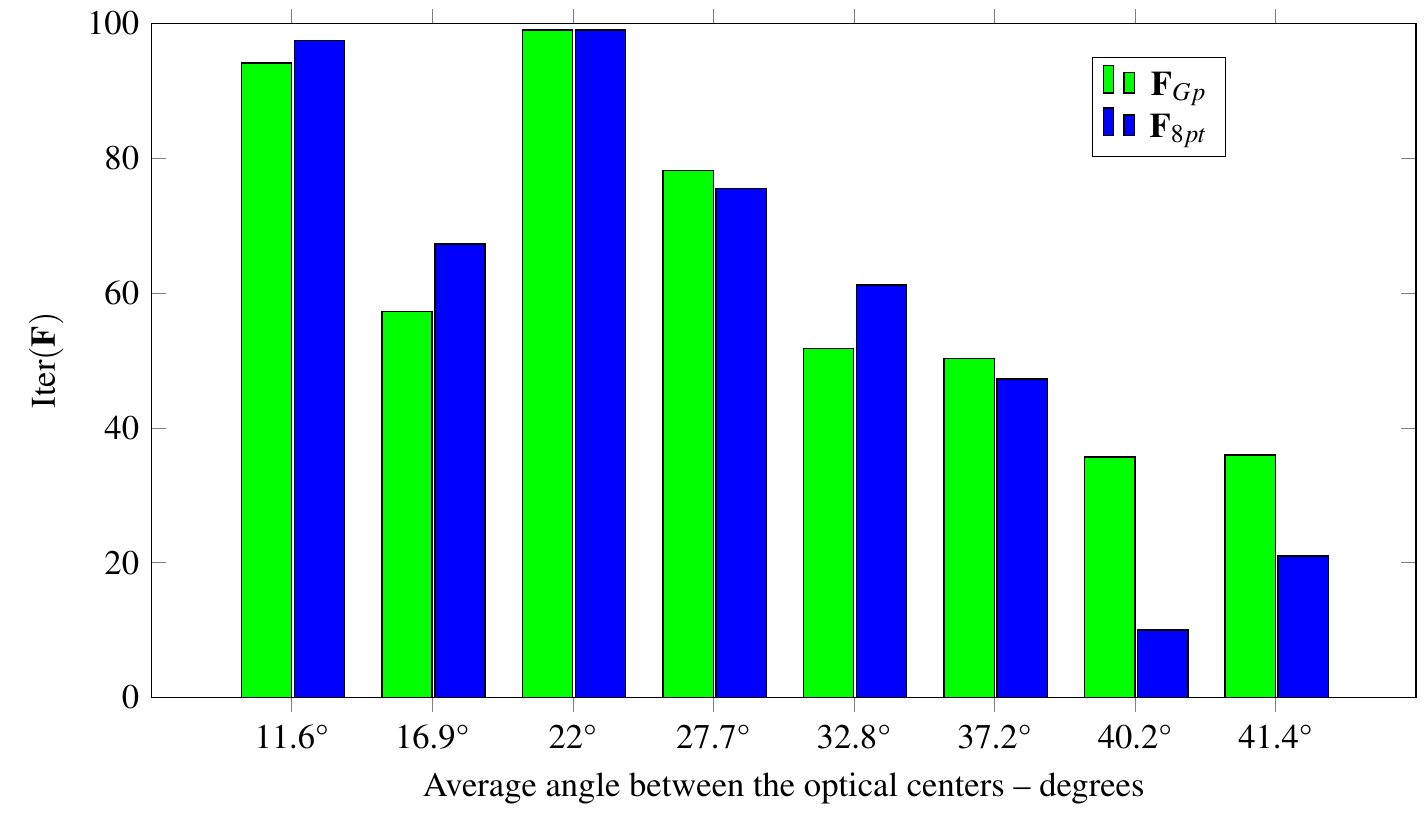}
\caption{\label{figure17} Number of iterations performed by Bundle-Adjustment to converge measured against movement amplitude for the \textit{House} series}
\end{figure*}

For the set of $36$ images of the \emph{Dinosaur} series and $9$ images of the {\em house} series, we tested the influence of motion amplitude between a pair of image on the quality of the resulting estimates obtained by both methods. For this purpose, we had both estimates with all possible motions $((0,1), (1,2), (2,3),\ldots)$ with $1$-image distance, then all possible motions $((0,2), (1,3),(2,4),\ldots)$ with $2$-image distance, and so on. With this process, we can measure the influence of the average angle on the quality of the fundamental matrix estimated by both methods. Figures from \ref{figure12} to \ref{figure17} shows the average of re-projection errors and the average of number of iterations with respect to average angle for the two series. The re-projection error after Bundle Adjustment is always smaller with the global method and with always a smaller number of iterations. Next, the larger the movement the more the solution by both methods deteriorates. But the 
deterioration is larger for the 8-point method than for the global method. One may also observe that in some cases the re-projection error before Bundle Adjustment is in favor of the 8-point method. In analogy with the real cases studied before, this may be due to the fact that for these cases the initial guess $\Mat{F}_{Gp}$ is in a basin of attraction with a better local minimum than in the basin of attraction associated with $\Mat{F}_{8Pt}$, but the `distance' between the initial guess and the corresponding local minimizer is larger for $\Mat{F}_{Gp}$ than for $\Mat{F}_{8Pt}$. Indeed in such cases the number of iterations is larger for $\Mat{F}_{Gp}$ than for $\Mat{F}_{8Pt}$.

\section{Conclusion}
We have studied the problem of estimating globally the fundamental matrix over nine parameters and under  rank and normalisation constraints. We have proposed a polynomial-based approach which 
enables one to estimate the fundamental matrix with good precision. More generally, we have shown how to modify the constraints on the numerical certificate of optimality to obtain fast and robust convergence. The method converges in a reasonable amount of time compared to other global optimization methods. 

From computational experiments conducted on both simulated and real data we conclude that the global method always provides an accurate initial estimation for the subsequent bundle adjustment step. Moreover,  we have shown that if the eight-point method has a  lower computational cost, its resulting estimate is frequently far from the global optimum obtained by the global method.



\begin{thebibliography}{xx}

\bibitem{longuet-higgins:81a}
H.~C. Longuet-Higgins.
\newblock {A computer algorithm for reconstructing a scene from two
  projections}.
\newblock {\em Nature}, 293:133--135, Sep 1981.

\bibitem{salvi:03}
X.~Armangu\'e and J.~Salvi.
\newblock Overall view regarding fundamental matrix estimation.
\newblock {\em Image and Vision Computing}, 21:205--220, 2003.

\bibitem{tsai:84a}
R.~Y. Tsai and T.~S. Huang.
\newblock {Uniqueness and estimation of three-dimensional motion parameters of
  rigid objects with curved surfaces}.
\newblock {\em IEEE Transaction on Pattern Analysis and Machine Intelligence},
  6:13--26, 1984.

\bibitem{hartley:95a}
R.~Hartley.
\newblock {In Defence of the 8-point Algorithm}.
\newblock In {\em 5$^{th}$ International Conference on Computer Vision
  (ICCV'95)}, pages 1064--1070, Boston (MA, USA), Jun 1995.

\bibitem{Torr02}
P.~H.~S. Torr and A.~W. Fitzgibbon.
\newblock Invariant fitting of two view geometry or "in defiance of the 8 point
  algorithm".
\newblock 2002.

\bibitem{Chojnacki03}
W.~Chojnacki, M.~J. Brooks, and A.~Van Den~Hengel.
\newblock Revisiting hartley's normalized eight-point algorithm.
\newblock {\em IEEE Transactions on Pattern Analysis and Machine Intelligence},
  25:1172--1177, 2003.

\bibitem{BartoliSturm:04}
A.~Bartoli and P.~Sturm.
\newblock Nonlinear estimation of the fundamental matrix with minimal
  parameters.
\newblock {\em IEEE Trans. Pattern Anal. Mach. Intell.}, 26(3):426--432, March
  2004.

\bibitem{Lasserre:01}
J.~B. Lasserre.
\newblock Global optimization with polynomials and the problem of moments.
\newblock {\em SIAM Journal on Optimization}, 11:796--817, 2001.

\bibitem{gloptipoly:02}
D.~Henrion and J.B. Lasserre.
\newblock Gloptipoly: Global optimization over polynomials with {Matlab} and
  {SeDuMi}.
\newblock {\em Proceedings of IEEE Conference on Decision and Control}, Dec.
  2002.

\bibitem{Xiao:10}
X.~Xuelian.
\newblock New fundamental matrix estimation method using global optimization.
\newblock In {\em Proceedings of International Conference on Computer
  Application and System Modeling (ICCASM)}, Taiyuan, China, 22--24 October
  2010.

\bibitem{Chesi:02}
G.~Chesi, A.~Garulli, A.~Vicino, and R.~Cipolla.
\newblock Estimating the fundamental matrix via constrained least-squares: A
  convex approach.
\newblock {\em IEEE Transactions on Pattern Analysis and Machine Intelligence},
  24:397--401, 2002.

\bibitem{KahlHenrion:07}
F.~Kahl and D.~Henrion.
\newblock Globally optimal estimates for geometric reconstruction problems.
\newblock {\em International Journal of Computer Vision}, 74:3--15, 2007.

\bibitem{salviPhd:99}
J.~Salvi.
\newblock {\em An Approach to Coded Structured Light to Obtain Three
  Dimensional Information}.
\newblock PhD thesis, University of Girona, 1999.

\bibitem{PeiChen:10}
P.~Chen.
\newblock Why not use the lm method for fundamental matrix estimation?
\newblock {\em IET computer vision}, 4(4):286--294, 2010.

\bibitem{luong:96a}
Q.-T. Luong and O.~Faugeras.
\newblock {The fundamental matrix: Theory, algorithms, and stability analysis}.
\newblock {\em International Journal of Computer Vision}, 17(1):43--76, 1996.

\bibitem{Zhang:98}
Z.~Zhang.
\newblock Determining the epipolar geometry and its uncertainty : A review.
\newblock {\em International Journal of Computer Vision}, 27:161--195, 1998.

\bibitem{Hengel02anew}
A.~Van Den~Hengel, W.~Chojnacki, M.~J. Brooks, and D.~Gawley.
\newblock A new constrained parameter estimator: Experiments in fundamental
  matrix computation.
\newblock In {\em In Proceedings of the 13th British Machine Vision
  Conference}, September 2002.

\bibitem{Chojnacki02anew}
W.~Chojnacki, M.~J. Brooks, A.~Van Den~Hengel, and D.~Gawley.
\newblock A new approach to constrained parameter estimation applicable to some
  computer vision problems.
\newblock In {\em Statistical Methods in Video Processing Workshop held in
  conjunction with ECCV'02}, Copenhagen, Denmark, 2002.

\bibitem{Hartley-Zisserman:03}
R.~Hartley and A.~Zisserman.
\newblock {\em Multiple View Geometry in Computer Vision}.
\newblock Cambridge University Press, 2003.

\bibitem{Torr97thedevelopment}
P.~H.~S. Torr and D.~W. Murray.
\newblock The development and comparison of robust methods for estimating the
  fundamental matrix.
\newblock {\em International Journal of Computer Vision}, 24:271--300, 1997.

\bibitem{Glover:86}
F.~Glover.
\newblock Future paths for integer programming and links to artificial
  intelligence.
\newblock {\em Comput. Oper. Res.}, 13(5):533--549, 1986.

\bibitem{Jaumard:90}
P.~Hansen and B.~Jaumard.
\newblock Algorithms for the maximum satisfiability problem.
\newblock {\em Computing}, 44(4):279--303, 1990.

\bibitem{Lang:60}
A.~H. Land and A.~G. Doig.
\newblock An automatic method of solving discrete programming problems.
\newblock {\em Econometrica}, 28:497--520, 1960.

\bibitem{hansen:03}
E.~R. Hansen and G.~W. Walster.
\newblock {\em Global optimization using interval analysis, 2nd Edn.}
\newblock Marcel Dekker, 2003.

\bibitem{Kearfott:09}
R.~E. Moore, R.~B. Kearfott, and M.~J. Cloud.
\newblock {\em Introduction to Interval Analysis}.
\newblock SIAM, 2009.

\bibitem{Weise:07}
T.~Weise.
\newblock {\em Global Optimization Algorithms - Theory and Application}.
\newblock Thomas Weise, 2007-05-01 edition, 2007.

\bibitem{kahl:05}
F.~Kahl and D.~Henrion.
\newblock Globally optimal estimates for geometric reconstruction problems.
\newblock In {\em Proc. IEEE Int. Conf. Computer Vision (ECCV)}, 2005.

\bibitem{Zheng:13}
Y.~Zheng, S.~Sugimoto, and M.~Okutomi.
\newblock A practical rank-constrained eight-point algorithm for fundamental
  matrix estimation.
\newblock In {\em Computer Vision and Pattern Recognition, 2013. CVPR 2013.
  IEEE Conference on}, June 2013.

\bibitem{lasserre00}
J.~B. Lasserre.
\newblock Optimisation globale et th\'eorie des moments.
\newblock {\em Comptes Rendus de l'Acad\'emie des Sciences}, 331:929--934,
  2000.

\bibitem{LasserreBook:10}
J.~B. Lasserre.
\newblock {\em Moments, positive polynomials and their applications}, volume~1
  of {\em Imperial College Press Optimization Series}.
\newblock Imperial College Press, 2010.

\bibitem{Laurent:08}
M.~Laurent.
\newblock Sums of squares, moment matrices and optimization over polynomials.
  {I}n {E}merging {A}pplications of {A}lgebraic {G}eometry.
\newblock {\em IMA Volumes in Mathematics and its Applications}, 149:157--270,
  2009.

\bibitem{Borchers99csdp}
B.~Borchers.
\newblock Csdp, a c library for semidefinite programming.
\newblock {\em Optimization Methods and Software}, 11:613--623, 1999.

\bibitem{dsdp-user-guide}
S.~J. Benson and Y.~Ye.
\newblock {DSDP5} user guide --- software for semidefinite programming.
\newblock Technical Report ANL/MCS-TM-277, Mathematics and Computer Science
  Division, Argonne National Laboratory, Argonne, IL, September 2005.
\newblock \url{http://www.mcs.anl.gov/~benson/dsdp}.

\bibitem{Fujisawa95sdpa}
K.~Fujisawa, M.~Fukuda, M.~Kojima, K.~Nakata, M.~Nakata, M.~Yamashita,
  K.~Fujisawa, M.~Fukuda, and K.~Kobayashi.
\newblock Sdpa (semidefinite programming algorithm) -- user's manual.
\newblock Technical report, Dept. Math. and Comp. Sciences -- Tokyo Institute
  of Technology, 2003.

\bibitem{sdpt3UserGuide}
K.C. Toh, M.J. Todd, and R.H. Tutuncu.
\newblock Solving semidefinite-quadratic-linear programs using sdpt3.
\newblock {\em Mathematical Programming}, 95:189--217, 2003.

\bibitem{S98guide}
J.F. Sturm.
\newblock Using {SeDuMi} 1.02, a {MATLAB} toolbox for optimization over
  symmetric cones.
\newblock {\em Optimization Methods and Software}, 11--12:625--653, 1999.
\newblock Version 1.05 available from {\texttt{http://fewcal.kub.nl/sturm}}.

\bibitem{Nie2012}
J.~Nie.
\newblock Optimality conditions and finite convergence of {L}asserre's
  hierarchy.
\newblock Technical report, Dept. Math. -- Univ. of California at San Diego,
  2012.

\bibitem{gloptipoly03}
D.~Henrion, J.B Lasserre, and J.~L\"ofberg.
\newblock Gloptipoly 3: moments, optimization and semidefinite programming.
\newblock {\em Optim. Methods and Software}, 24:761--779, 2009.

\bibitem{YALMIP}
J.~Löfberg.
\newblock Yalmip : A toolbox for modeling and optimization in {MATLAB}.
\newblock In {\em Proceedings of the CACSD Conference}, Taipei, Taiwan, 2004.

\bibitem{sdpt3}
K.C. Toh, M.J. Todd, and R.H. Tutuncu.
\newblock Sdpt3 --- a matlab software package for semidefinite programming.
\newblock {\em Optimization Methods and Software}, 11:545--581, 1999.

\bibitem{Luong:1996}
Q.-T. Luong and T.~Vi\'{e}ville.
\newblock Canonical representations for the geometries of multiple projective
  views.
\newblock {\em Comput. Vis. Image Underst.}, 64:193--229, September 1996.

\bibitem{hartleysturmcviu}
R.~Hartley and P.~Sturm.
\newblock Triangulation.
\newblock {\em Computer Vision and Image Understanding}, 68(2):146--157, 1997.

\end{thebibliography}
\end{document}